\documentclass[12pt]{article} 

\usepackage{url}
\usepackage{a4wide}
\usepackage{amsmath}
\usepackage{amssymb}
\usepackage{color}
\usepackage{bm}
\usepackage{graphicx}
\usepackage{hyperref}

\let\pa\partial

\newcommand{\R}{\mathbb{R}}

\renewcommand\Re{\operatorname{Re}}

\newcommand{\diag}{\operatorname{diag}}

\newcommand{\Ekininc}{E_{\mathrm{kin}}^{\mathrm{inc}}}

\author{
	Jan-Frederik Mennemann
  	\footnote{jfmennemann@gmx.de,
	Institute for Analysis and Scientific Computing, Vienna University of Technology, 
	Wiedner Hauptstra\ss e 8--10, 1040 Vienna, Austria
	}
  	\and
  	Ansgar J\"ungel
  	\footnote{juengel@tuwien.ac.at,
	Institute for Analysis and Scientific Computing, Vienna University of Technology, Wiedner 		Hauptstra\ss e 8--10, 1040 Vienna, Austria
	}
}

\title{Perfectly Matched Layers versus discrete transparent boundary conditions in quantum device simulations}

\begin{document}

\maketitle


\begin{abstract} 
Discrete transparent boundary conditions (DTBC) and the 
Perfectly Matched Layers (PML) method for the realization of open boundary conditions
in quantum device simulations are compared, based on the stationary and 
time-dependent Schr\"odinger equation. The comparison includes scattering
state, wave packet, and transient scattering state simulations in one and
two space dimensions. 
The Schr\"o\-dinger equation is discretized by a second-order Crank-Nicolson
method in case of DTBC. For the discretization with PML,
symmetric second-, fourth, and sixth-order spatial approximations as well as
Crank-Nicolson and classical Runge-Kutta time-integration methods
are employed. In two space dimensions, a ring-shaped quantum waveguide device is
simulated in the stationary and transient regime. 
As an application, a simulation of the Aharonov-Bohm effect in this device is performed,
showing the excitation of bound states localized in the ring region.
The numerical experiments show that the results obtained from PML are 
comparable to those obtained using DTBC, while keeping the high numerical 
efficiency and flexibility as well as the ease of implementation of the former method.
\end{abstract}


\medskip


\section{Introduction}

The electron transport in nanoscale electronic devices, 
in which inelastic collisional effects may be neglected, 
can be modeled by the stationary or  time-dependent Schr\"odinger equation.
The semiconductor device is typically connected to semi-infinite leads
which describe the electric contacts. The aim is to solve the Schr\"odinger
problem in the bounded device domain instead in the whole space. This
makes it necessary to prescribe appropriate open boundary conditions at the interface 
between the leads and the active region of the device to avoid unphysical
reflections at the boundary. In the literature, several methods have been
proposed to derive transparent boundary conditions. Analytical transparent
boundary conditions are nonlocal in time, and their numerical implementation
requires some care (see the review \cite{AABES08} and references therein).
Moreover, inadequate discretizations may introduce strong reflections at the
boundary. 

In this paper, we implement and compare two approaches: the discrete
transparent boundary conditions (DTBC), which are completely reflection-free
at the boundary but nonlocal in time, and 
the Perfectly Matched Layers (PML), which involve an artificial boundary layer
but which can be implemented very efficiently.

In the context of finite-difference discretizations, 
DTBC were derived by Arnold \cite{Arn98}.
They yield unconditionally stable reflection-free numerical discretizations.
DTBC include the discrete convolution of the unknown function with a given 
kernel, whose numerical computation is rather involved. The evaluation
can be significantly accelerated by approximating the kernel by a finite sum
of exponentials which decay with respect to time \cite{AES03}. 
The limit of
vanishing spatial approximation parameters in the DTBC coincides with the
temporally semi-discrete transparent boundary conditions of \cite{ScDe95,LuSc02}.
DTBC for the Schr\"odinger equation were also derived for finite-element
\cite{ZlZl12} and splitting higher-order schemes \cite{DZR13}.

The PML approach was introduced by B\'erenger \cite{Ber94} for absorbing
boundaries for wave equations. The idea is to replace the absorbing 
boundary condition by an absorbing boundary layer which damps out waves 
using a damping function. The problem of having reflections from the
absorber boundary was handled by B\'erenger by constructing a special
absorbing medium. PML can be seen as the result of a complex coordinate
transformation, being essentially a continuation of the Schr\"odinger
equation into complex spatial coordinates \cite{ChWe94}.
Later, PML were derived and analyzed for many other
equations, like wave and Helmholtz equations \cite{TuYe98} and
Schr\"odinger equations \cite{Col97}. Nonlinear Schr\"odinger problems
were investigated in, for instance, \cite{CLLML07,Zhe07,ZhCh10}. 
The simulations include (stationary) scattering state calculations \cite{CLLML07} or
wave packets \cite{Zhe07,ZhCh10}.
However, transient scattering states which are needed to describe the time-dependent
behavior of quantum devices are not considered.

In a transient scattering state simulation, the initial wave function is given by a 
(stationary) scattering state.
When the transient simulation is started, the external potential is allowed to 
change over time and hence, the scattering state starts to evolve in time too.

Our aim is to provide a careful numerical comparison between DTBC and PML
in the finite-difference context. 
We show how PML can be applied in a variety of physical situations:
(stationary) scattering states, wave packets, incoming waves, and transient 
scatterings states.
In the following, we describe our approach and the main results in more detail.

Since the open boundary problem for quantum waveguides in several space
dimensions can be reduced to the one-dimensional case, we consider one-dimensional
simulations first. The Schr\"odinger equation is discretized using symmetric
finite differences of second, fourth, and sixth order.
The numerical solutions obtained from (second-order) DTBC and (second- and
higher-order) PML are compared in a series of simulations. 
It turns out
that for second-order discretizations, PML compete well with DTBC for 
medium to large energies ($10\ldots1000$\,meV).
For very small energies ($10^{-3}\ldots 1$\,meV), 
second-order DTBC perform significantly better than second-order
PML, but the numerical errors for higher-order PML schemes are comparable
to that for second-order DTBC at small energies and they are much smaller
for medium to large energies.

Then we turn to two-dimensional quantum waveguide simulations.
The stationary scattering state problem has been solved in \cite{LeKi90} 
using linear finite elements along with exact transparent boundary conditions.
The transient case is considered in \cite{ArSc08,ScAr08}
using DTBC based on the Crank-Nicolson scheme.
However, transient scattering states are not considered and the cross sections of the leads need to be infinite square well potentials.
In this paper, we show how to remove these limitations.
Although we consider only ring-shaped two-terminal quantum waveguides,
our approach can also be applied to more complicated multi-terminal devices.

The implementation of open boundary conditions using PML works analogously
to the one-dimensional case. 
Contrary to DTBC, a decomposition into cross-sectional waveguide eigenstates 
is not required, which simplifies the implementation significantly. 
Our simulations indicate that the numerical error which results from the PML 
is of the same
order as that resulting from the approximation of the Schr\"odinger equation
itself. In contrast to DTBC, which are tailored specifically to numerical
methods, PML can be applied in a more flexible way. As an example, we employ
second-, fourth-, and sixth-order finite-difference formulas to approximate the
spatial derivatives. For the time integration, we use the Crank-Nicolson 
approximation and classical Runge-Kutta methods. 

The Runge-Kutta approach has the drawback that the resulting spatio-temporal
discretization is only conditionally stable and that the mass of the particles
is not conserved exactly. Our numerical results suggest, however, that these
issues may be overestimated. The advantages are the competitive computing times,
the high accuracy for higher-order schemes, and the easy implementation.
The combination of the classical Runge-Kutta method,
higher-order spatial discretizations, and PML appears to be a very efficient
approach for transient quantum device simulations.

Finally, we apply the above mentioned methods in a simulation of the 
Aharonov-Bohm effect \cite{AhBo59}.
Transient simulations of the Aharonov-Bohm effect were also studied in the literature.
For instance, the electron transport in a quantum ring, using an expansion of the wave function in a basis of Gaussians and a finite-difference approach, was considered
in \cite{SzPe05}.
Simulations of ring-quantum interference transistors, 
where the Aharonov-Bohm effect is studied in dependence of externally applied
electrostatic potentials, are shown in \cite{HeJa00}.
However, these works do not include transient scattering state simulations.
As far as we know, we present in this paper the first two-dimensional transient 
scattering state simulations of the Aharonov-Bohm effect in a ring-shaped quantum 
waveguide. Remarkably, we find that fast variations of the external magnetic field lead to the excitation of bound states which are localized in the ring region.

The paper is organized as follows. Section \ref{sec:one_dimensional_simulations} 
covers stationary
scattering states, wave packets, and transient scattering states in one
space dimension. We recall the discretization of the Schr\"odinger equation
using DTBC, detail the approximation employing PML, 
and compare both methods numerically.
In Section \ref{sec:quantum_waveguide_simulations}, 
we consider two-dimensional ring-shaped quantum waveguide devices.
The derivation of DTBC and PML for the stationary and transient problems is detailed.
Section \ref{sec.AB} is devoted to the Aharonov-Bohm effect in the
stationary and transient regime.


\section{One-dimensional simulations}
\label{sec:one_dimensional_simulations}

\subsection{Scattering states}
\label{subsec:scattering_states_1d}

Scattering states in one-dimensional simulations represent a beam of electrons
injected at the left or right lead of the device.
This beam of particles is identified with a wave function which solves the
stationary Sch\"odinger equation
\begin{equation}\label{eq:stationary_schroedinger_equation_1d}
  -\frac{\hbar^2}{2 m^\star} \frac{d^2}{dx^2}\phi(x) + V(x) \phi(x) = E \phi(x), 
	\quad x \in \R,
\end{equation}
subject to the boundary condition at infinity 
that the incoming wave function is a plane wave. 
Here, $E$ denotes the total energy, $V$ the potential energy, $\hbar$ the
reduced Planck constant, and $m^\star$ the effective mass of the electrons 
in the semiconductor.
In all subsequent simulations, we choose $m^\star=0.067\,m_e$, corresponding to the 
effective mass of electrons in GaAs, with $m_e$ being the electron mass at rest.

The device in $(0,L)$ is assumed to be connected to the semi-infinite leads
$(-\infty,0]$ and $[L,\infty)$. In the exterior domain $(-\infty,0]\cup[L,\infty)$,
the potential energy is assumed to be constant, i.e.\
$V(x)=V_\ell$ for $x \leq 0$ and $V(x)=V_r$ for $x \geq L$.
Without loss of generality we set $V_\ell=0$ and $V_r = -e U$,
where $e$ denotes the elementary charge and $U$ is the applied voltage 
at the right contact.

The energy of an electron injected at the left contact is given by 
$E(k) = \hbar^2 k^2 / (2 m^\star)$, where $k>0$ denotes the electron 
wave number in the left lead. 
The incoming electron is represented by a plane wave $\exp(i k x)$ 
traveling to the right.
Accordingly, transparent boundary conditions read as \cite{Arn01}
\begin{equation}\label{eq:stationary_transparent_boundary_conditions_1d}
  \phi'(0) + i k \phi(0) = 2 i k, \quad 
  \phi'(L) = i \sqrt{2 m^\star (E -V_r) / \hbar^2} \phi(L).
\end{equation}
Electrons injected at the right contact traveling to the left are treated analogously.

\medskip\noindent
{\bf Discrete transparent boundary conditions.}
A symmetric second-order finite-diffe\-rence discretization 
of the stationary Schr\"odinger equation 
\eqref{eq:stationary_schroedinger_equation_1d} on the equidistant grid 
$x_j=j\triangle x$, $j \in \mathbb{Z}$, with $x_J=L$ and $\triangle x > 0$
is given by
\begin{equation}\label{eq:discrete_stationary_schroedinger_equation_1d}
  -\frac{\hbar^2}{2 m^\star} 
  \frac{\phi_{j-1} - 2 \phi_{j} + \phi_{j+1}}{(\triangle x)^2}
  + V_j \phi_j = E \phi_j.
\end{equation}
As before, the potential energy in the semi-infinite leads is assumed to be constant,  
i.e., $V_j=0$ for $j \leq 0$ and $V_j=-eU$ for $j \geq J$.
We seek for a solution of \eqref{eq:discrete_stationary_schroedinger_equation_1d}
restricted to the grid
\begin{equation}\label{eq:X_dtbc}
  X_{\mathrm{DTBC}} = \left\{x_j = j \triangle x, \; j = 0,\ldots,J \right\},
\end{equation}
and hence we have to specify boundary conditions at $x_0=0$ and $x_J=L$.
It is well known that a direct discretization of the transparent boundary conditions 
in \eqref{eq:stationary_transparent_boundary_conditions_1d} via 
standard finite-difference formulas may lead to spurious oscillations 
in the numerical solution.
These oscillations may be eliminated using DTBC \cite{Arn01}.
The idea is to derive the boundary conditions on the discrete level of equation
\eqref{eq:discrete_stationary_schroedinger_equation_1d}.
In the semi-infinite leads, the potential energy is assumed to be constant and thus,
\eqref{eq:discrete_stationary_schroedinger_equation_1d} 
admits two solutions of the form
$\phi_j = \alpha_{\ell,r}^j$,
where the numbers
$$
  \alpha_{\ell,r} = 1 - \frac{m^\star E_{\ell,r}^{\mathrm{kin}}(k) 
	(\triangle x)^2}{\hbar^2} 
  \pm i 
  \sqrt{\frac{2 m^\star E_{\ell,r}^{\mathrm{kin}}(k)(\triangle x)^2}{\hbar^2}
	-\frac{(m^\star)^2E_{\ell,r}^{\mathrm{kin}}(k)^2 (\triangle x)^4}{\hbar^4}}
$$
depend on the kinetic energy $E_{\ell,r}^{\mathrm{kin}}(k) = E(k)-V_{\ell,r}$ 
in the left or right lead, respectively.

Depending on $|\alpha_{\ell,r}| \lessgtr 1$, the solutions are exponentially 
decreasing or increasing.
In case $|\alpha_{\ell,r}| = 1$ the solutions are discrete plane waves (see below)
unless $\alpha_{\ell,r} = 1$ in which case they are constant.
For a wave function injected at the left contact traveling to the right, 
the solution  to \eqref{eq:discrete_stationary_schroedinger_equation_1d}
is a superposition of an incoming and a reflected discrete plane wave,
$\phi_j = A \alpha_\ell^j + B \alpha_\ell^{-j}$ in the left contact,
and a transmitted wave $\phi_j = C \alpha_r^j$ in the right contact.
The amplitude of the incoming discrete plane wave in the left lead is set to $A=1$.
Eliminating $B$ and $C$ yields the desired DTBC
\begin{equation}
  \label{eq:dtbc_stationary_1d}
  \phi_0 -\alpha_\ell \phi_1 = 1-\alpha_\ell^2, \quad 
  \alpha_r \phi_{J-1} - \phi_J = 0.
\end{equation}
Here, we implicitly assumed $V_j=0$ for $j\leq 1$ and $V_j=-eU$ for $j\geq J-1$
but we could just as easily state the boundary conditions in terms of 
$\phi_{-1}, \phi_0$ and $\phi_J, \phi_{J+1}$.

On the discrete level of equation 
\eqref{eq:discrete_stationary_schroedinger_equation_1d},
we need to replace the continuous $E$--$k$--relation
\begin{equation}
  \label{eq:continuous_E_k_relation}
  k = \sqrt{2 m^\star E^\mathrm{kin}} / \hbar
\end{equation}  
by the discrete $E$--$k$--relation
\begin{equation}
  \label{eq:discrete_E_k_relation}
  k = \arccos \left(1 - m^\star (\triangle x / \hbar)^2 E^\mathrm{kin} 
	\right) / \triangle x,
\end{equation}
which approximates the continuous relation up to second order.
Consequently, the discrete solutions can be written as
$\phi_j = \exp(\pm i k j \triangle x)$ and hence for $|\alpha_{0,J}^\pm| = 1$, 
they are called discrete plane waves.
Plane waves injected at the right contact traveling to the left are treated
analogously. For further details we refer to \cite{Arn01}.

\medskip\noindent
{\bf Perfectly Matched Layers.}
PML can be formulated using the complex coordinate transformation
\begin{equation}
  \label{eq:complex_coordinate_transformation}
  x \mapsto x + e^{i \pi / 4} \int^x \sigma(x')\; dx',
\end{equation}
with the absorption function
$$
  \sigma(x) = 
  \begin{cases}
	\begin{aligned}
  \sigma_0 (x^\star_{\ell}-x)^p & \quad\mbox{for } x < x_\ell^\star, \\
  0 & \quad\mbox{for } x_\ell^\star \leq x \leq x_r^\star, \\
  \sigma_0 (x-x^\star_r)^p & \quad\mbox{for } x_r^\star < x.
	\end{aligned}
  \end{cases}
$$
This coordinate transformation corresponds to the substitution
\begin{equation}
  \label{eq:pml_complex_function_c}
  \frac{d}{dx} \rightarrow c(x) \frac{d}{dx}, \quad 
	c(x):= \frac{1}{1+e^{i \pi / 4} \sigma(x)},
\end{equation}
such that the stationary Schr\"odinger-PML equation reads as
\begin{equation}
  \label{eq:stationary_schroedinger_pml_equation_1d}
  -\frac{\hbar^2}{2 m^\star} c(x) \frac{d}{d x} 
	\left( c(x) \frac{d}{d x} \phi(x) \right) + V(x) \phi(x) = E \phi(x),
	\quad x\in{\mathbb R}.
\end{equation}
This agrees with the original stationary Schr\"odinger equation 
\eqref{eq:stationary_schroedinger_equation_1d}
for $x \in [x_\ell^\star, x_r^\star]$, and hence we require that
$[0,L] \subset [x_\ell^\star, x_r^\star]$.
Outside of $[x_\ell^\star, x_r^\star]$, propagating waves are damped 
exponentially fast,
\begin{equation}
  \label{eq:propagating_wave_in_pml}
  \exp(i k x) \rightarrow \exp(i k x)
  \exp \left( i k e^{i\pi/4} \int^x \sigma(x') \; dx' \right),
\end{equation}
with their distance to the points $x_{\ell,r}^\star$.
Denoting the thickness of the PML by $d_0$, 
the computational domain is given by 
$[x_\ell^\star-d_0, x_r^\star+d_0]$.
We use cubic absorption functions ($p=3$), which we found to give slightly better 
results than the quadratic functions used in \cite{Zhe07}.

A wave propagating through a PML is expected to be practically zero when it hits 
the boundary of the computational domain,
provided $d_0$ is sufficiently large.
Therefore, either Dirichlet or Neumann boundary conditions can be imposed 
at the boundary points $x^\star_\ell-d_0$ and $x^\star_r+d_0$.
We discuss this issue below in more detail.

Propagating waves with different wave numbers $k$ experience different 
attenuation according to \eqref{eq:propagating_wave_in_pml}.
A uniform attenuation independent of the wave number can be obtained by including
$1/k$ in the absorption function.
This strategy is typically used in stationary wave problems \cite{SiTu04}.
Since we are concerned mainly with transient simulations, 
we employ a different strategy.

Any wave packet can be thought of as a superposition of propagating waves 
of different energy or different wave number $k$.
For that reason, we seek for PML which are able to treat all incoming waves 
simultaneously. In this spirit, the factor $\exp(i\pi/4)$ appearing in 
\eqref{eq:complex_coordinate_transformation}
is meant to give a good effect on average as explained in \cite{Zhe07}.
The thickness of the PML will be 
comparatively large as they are not optimized to absorb waves of a single energy.
In the simulations below, we fix the thickness of the PML first.
Then we choose the absorption strength factor $\sigma_0$ 
such that the numerical error becomes sufficiently small 
for a given range of energies.

We use symmetric finite-difference formulas to approximate the modified 
spatial derivative $c(x)\pa_x(c(x)\pa_x)=c(x)c'(x)\pa_x + c(x)^2\pa_x^2$
in \eqref{eq:stationary_schroedinger_pml_equation_1d}.
Approximations of second, fourth and sixth order read as
\begin{subequations}
  \label{eq:D_tilde_1d}
  \begin{align}
  \begin{split}
  \label{eq:D_tilde_1d_2nd}
  \tilde{D}_x^{2,\mathrm{2nd}} \phi_j
  := c(x_j) c'(x_j) D_x^{1,\mathrm{2nd}} \phi_j + c^2(x_j) D_x^{2,\mathrm{2nd}} \phi_j,
  \end{split}
  \\
  \begin{split}
  \label{eq:D_tilde_1d_4th}
  \tilde{D}_x^{2,\mathrm{4th}} \phi_j
  := c(x_j) c'(x_j) D_x^{1,\mathrm{4th}} \phi_j + c^2(x_j) D_x^{2,\mathrm{4th}} \phi_j,
  \end{split}
  \\
  \begin{split}
  \label{eq:D_tilde_1d_6th}
  \tilde{D}_x^{2,\mathrm{6th}} \phi_j
  := c(x_j) c'(x_j) D_x^{1,\mathrm{6th}} \phi_j + c^2(x_j) D_x^{2,\mathrm{6th}} \phi_j,
  \end{split}
  \end{align}
\end{subequations}
where the abbreviations
\begin{subequations}
  \label{eq:D_1_2nd_4th_6th}
  \begin{align}
  \begin{split}
  \label{eq:D_1_2nd}
  D_x^{1,\mathrm{2nd}} \phi_j
  &:=(-\phi_{j-1} + \phi_{j+1})/(2 \triangle x),
  \end{split}
  \\
  \begin{split}
  \label{eq:D_1_4th}
  D_x^{1,\mathrm{4th}} \phi_j
  &:= (\phi_{j-2} - 8 \phi_{j-1} + 8 \phi_{j+1} - \phi_{j+2})/(12 \triangle x),
  \end{split}
  \\
  \begin{split}
  \label{eq:D_1_6th}
  D_x^{1,\mathrm{6th}} \phi_j
  &:=
  (-\phi_{j-3} + 9 \phi_{j-2} - 45 \phi_{j-1} + 45 \phi_{j+1} - 9 \phi_{j+2} 
	+ \phi_{j+3})/(60 \triangle x), 
  \end{split}
	\\
  \begin{split}
  \label{eq:D_2_2nd}
  D_x^{2,\mathrm{2nd}} \phi_j
  &:= (\phi_{j-1}-2 \phi_j + \phi_{j+1})/(\triangle x)^2,
  \end{split}
  \\
  \begin{split}
  \label{eq:D_2_4th}
  D_x^{2,\mathrm{4th}} \phi_j
  &:= (-\phi_{j-2} + 16 \phi_{j-1} - 30 \phi_{j} + 16 \phi_{j+1} 
	- \phi_{j+2})/(12 (  \triangle x)^2),
  \end{split}
  \\
  \begin{split}
  \label{eq:D_2_6th}
  D_x^{2,\mathrm{6th}} \phi_j
  &:= (2 \phi_{j-3} - 27 \phi_{j-2} + 270 \phi_{j-1} - 490 \phi_{j} + 270 \phi_{j+1}
  -27 \phi_{j+2} \\
  &\phantom{:=i}{}  + 2 \phi_{j+3}^{(n)})/(180 (\triangle x)^2),
  \end{split}
  \end{align}
\end{subequations}
are needed frequently in the following.
We further introduce the equidistant grid
\begin{equation}
\label{eq:X_pml}
X_{\mathrm{PML}}:=
\left\{
x_j = x_\ell^\star-d_0 + j \triangle x, \; j = 0,\ldots,J_\mathrm{PML}
\right\}
\end{equation}
with $x_{J_\mathrm{PML}} = x_r^\star+d_0$ and 
$X_\mathrm{DTBC} \subset X_\mathrm{PML}$.
Then the discrete stationary Schr\"odinger-PML equation reads
\begin{equation}
\label{eq:discrete_stationary_schroedinger_pml_equation_1d}
-\frac{\hbar^2}{2 m^\star} \tilde{D}_x^{2} \, \phi_j
+ V_j \phi_j = E \phi_j, 
\quad \tilde{D}_x^2 \in 
\{\tilde{D}_x^{2,\mathrm{2nd}}, \tilde{D}_x^{2,\mathrm{4th}}, \tilde{D}_x^{2,\mathrm{6th}}\}.
\end{equation}
At the boundaries of $X_\mathrm{PML}$, we impose homogeneous Dirichlet or Neumann 
boundary conditions.
In the latter case, $\tilde{D}_x^2$ is modified accordingly.

We still need to specify how to realize an incoming plane wave at, say, 
the left contact.
To simplify the presentation, we restrict ourselves to the second-order 
discretization $\tilde{D}_x^{2}=\tilde{D}_x^{2,\mathrm{2nd}}$.
Let $x_{j_0} = 0$ denote the boundary of the left contact.
The wave function in the left contact is given by the sum of an incoming and 
a reflected wave function $\phi_j = \phi_j^\mathrm{inc} + \phi_j^\mathrm{refl}$.
In the numerical implementation, we eliminate $\phi_j^\mathrm{inc}$ in the left lead.
Thus, to realize an incoming plane wave at $x_{j_0}$,
the finite difference equations for $j_0-1$ and ${j_0}$ need to be modified 
according to
\begin{subequations}
  \label{eq:incoming_wave_stationary_schroedinger_1d}
  \begin{align}
  -\frac{\hbar^2}{2 m^\star}
  \frac{\phi_{j_0-2} - 2 \phi_{j_0-1} 
	+ (\phi_{j_0}-\phi_{j_0}^\mathrm{inc})}{(\triangle x)^2} + V_{j_0-1} \phi_{j_0-1}
  &= E \phi_{j_0-1}, \\
  -\frac{\hbar^2}{2 m^\star}
  \frac{(\phi_{j_0-1}+\phi_{j_0-1}^\mathrm{inc}) 
	- 2 \phi_{j_0} + \phi_{j_0+1}}{(\triangle x)^2} + V_{j_0} \phi_{j_0}
  &= E \phi_{j_0}.
  \end{align}
\end{subequations}

We note that the absorption function of the left lead is only active for 
$x_j \leq x_\ell^\star$. Choosing $x_\ell^\star$ slightly smaller than $x_{j_0}$
ensures that the complex function $c$ does not show up in
\eqref{eq:incoming_wave_stationary_schroedinger_1d}.
Therefore, the final numerical problem becomes
\begin{equation*}
  -\frac{\hbar^2}{2 m^\star} \tilde{D}_x^{2,\mathrm{2nd}} \phi_j + (V_j - E)
	\phi_j = b_j,\quad j = 0,\ldots,J_\mathrm{PML},
\end{equation*}
where
$$
  b_j =
  \begin{cases}
  \begin{aligned}
  -\left( \hbar^2 / (2 m^\star (\triangle x)^2) \right) \phi_{j_0}^\mathrm{inc} 
	&\quad   \text{for } j = j_0-1, \\
  +\left( \hbar^2 / (2 m^\star (\triangle x)^2) \right) \phi_{j_0-1}^\mathrm{inc} 
	&\quad \text{for } j = j_0, \\
  0 &\quad \textrm{else}.
  \end{aligned}
  \end{cases}
$$
For higher-order approximations $\tilde{D}_x^2=\tilde{D}_x^{2,\mathrm{4th}}$ 
or $\tilde{D}_x^2=\tilde{D}_x^{2,\mathrm{6th}}$, we proceed in a similar way.
However, due to the extended finite-difference stencils, 
the vector $b$ involves four or six non-zero entries, respectively.

The incoming plane wave is given by $\phi_j^\mathrm{inc} = \exp(i k j \triangle x)$. 
In case of the second-order discretization, 
the wave number is related to the kinetic energy by the discrete 
$E$--$k$--relation \eqref{eq:discrete_E_k_relation}.
In case of the higher-order discretizations, we simply use 
\eqref{eq:continuous_E_k_relation}, since the corresponding discrete 
$k$--$E$--relations cannot easily be inverted and the differences 
would be small anyhow.

\medskip\noindent
{\bf Simulations.}
We consider the ramp-like potential energy
\begin{equation}
  \label{eq:ramp_like_potential}
  V(x) =
  \begin{cases}
  \begin{aligned}
  0 &\quad \text{for } x < a_0, \\
  -\dfrac{x-a_0}{a_1-a_0} e U &\quad \text{for } a_0 \leq x < a_1, \\
  -e U    &\quad \text{for } x \geq a_1,
  \end{aligned}
  \end{cases}
\end{equation}
with $a_0=40$\,nm, $a_1=80$\,nm, and $U=-25$\,mV.
The device extends from $0$\,nm to $L=120$\,nm.
The electrons are injected at the left contact traveling to the right.
At the left contact, the potential energy is zero, and hence the energy of the 
incoming electrons is given by the kinetic energy only, which is denoted by 
$\Ekininc$. We choose the width $d_0=40$\,nm of the PML, the absorption 
strength factor $\sigma_0=0.02$, and the mesh size $\triangle x=0.5$\,nm. 
Unless stated otherwise, the simulations in this subsection are performed using 
homogeneous Neumann boundary conditions at the end points of the PML.

\begin{figure}[thb]
	\centering
  	\includegraphics[width=150mm]{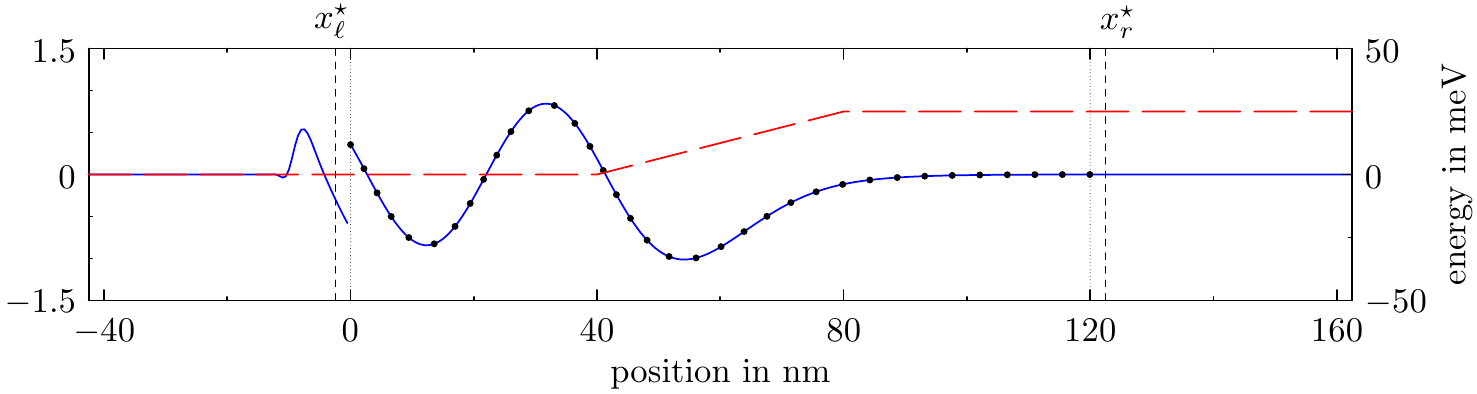}
  	\includegraphics[width=150mm]{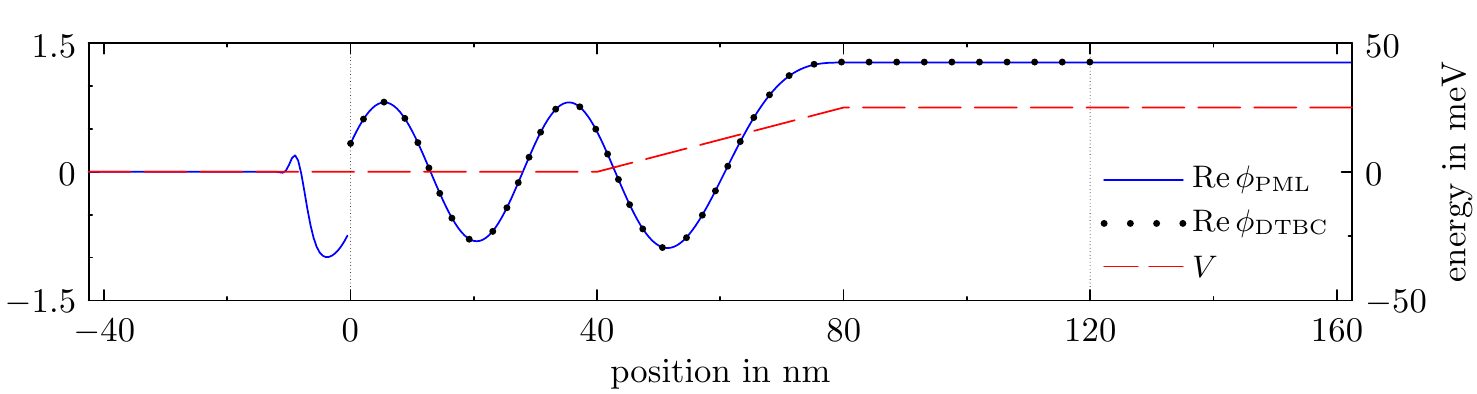}
  	\includegraphics[width=150mm]{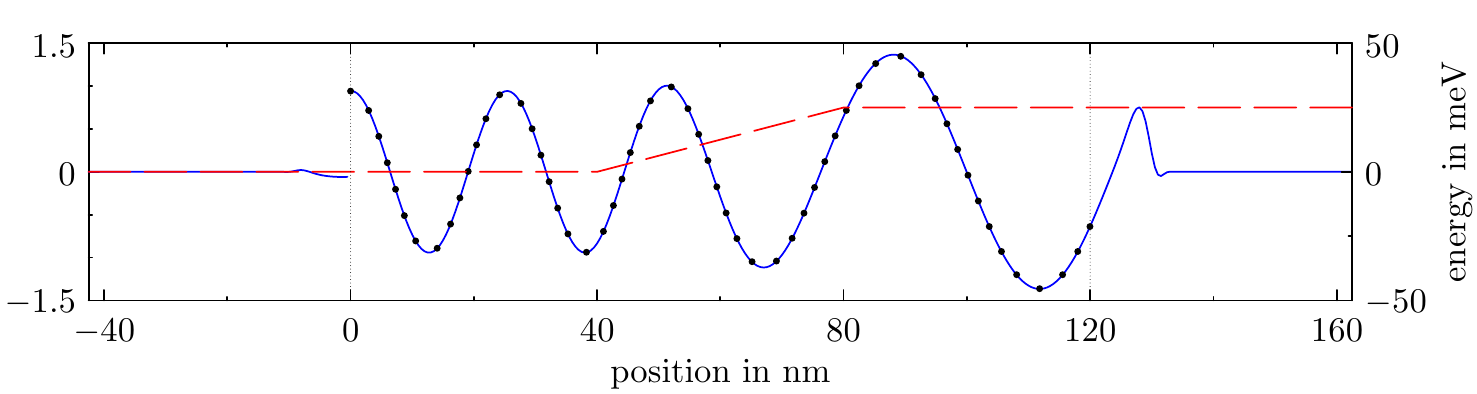}
  \caption{Scattering states for an electron injected at the left contact.
  The kinetic energy of the incoming electrons is
  $E_\mathrm{inc}^\mathrm{kin} = 15$\,meV (top), $25$\,meV (center), and
  $35$\,meV (bottom).
  The points $x_\ell^\star$ and $x_r^\star$ indicate the boundaries of the PML 
	in the left and right lead, respectively.}
\label{fig:scattering_states_1d}
\end{figure}

Figure \ref{fig:scattering_states_1d} shows the real parts 
of scattering states computed with DTBC (dotted line) and PML (solid line)
for different energies $\Ekininc$.
The distances of $x_\ell^\star$ and $x_r^\star$ to the device domain are 
a few times the mesh size $\triangle x$, which ensures that the discretization 
inside the device domain is not altered by the PML.
We note that the second-order discretizations with DTBC and PML 
coincide exactly inside the device domain. Indeed, 
Figure \ref{fig:scattering_states_1d} shows that the scattering states of both methods,
computed for different values of $\Ekininc$, can hardly by distinguished.

Figure \ref{fig:scattering_states_1d} (top) corresponds to the energy
$E_\mathrm{kin}^\mathrm{inc}=15$\,meV.
Since the potential energy at the right contact amounts to $25$\,meV, 
all incoming electrons are being reflected and then absorbed by the 
PML in the left contact.
The discontinuity in the wave function at $x=0$\,nm stems from the fact that 
the incoming plane wave in the left contact has been eliminated.

The middle figure corresponds to the limiting case $\Ekininc=25$\,meV, i.e.,
electrons which reach the contact on the right-hand side have zero kinetic energy.
Hence, the wave number becomes zero and according to 
\eqref{eq:propagating_wave_in_pml},
the PML in the right contact has no effect on the wave function.
Nonetheless, we obtain a reasonable approximation due to the 
Neumann boundary conditions imposed at the boundary of the computational domain.
However, the question arises whether the PML fails when
$E_\mathrm{inc}^\mathrm{kin}$ is very close to this critical energy.
This issue will be addressed below.

Figure \ref{fig:scattering_states_1d} (bottom)
corresponds to $E_\mathrm{kin}^\mathrm{inc}=35$\,meV.
Since the potential energy of the transmitted electrons is increased by $25$\,meV,
the kinetic energy is decreased by this amount which results in a
reduced wave number, or equivalently, in an increased wavelength.

\begin{figure}[t]
  \centering
  \includegraphics[width=75mm]{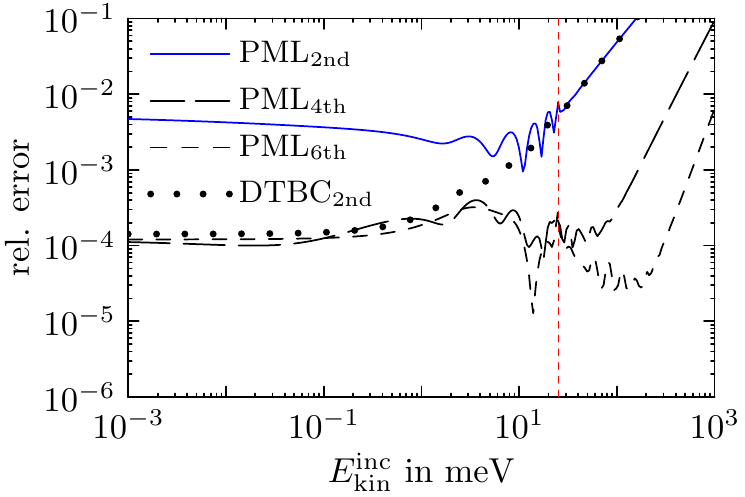}
  \hspace{2.5mm}
  \includegraphics[width=75mm]{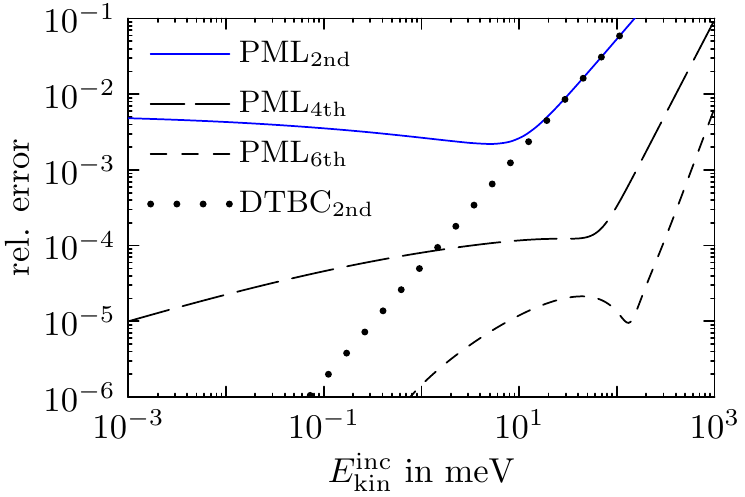}
  \caption{Relative errors of the scattering states as a function 
  of the kinetic energy with applied voltage $U=-25$\,mV (left) and $U=0$ (right).
  The critical energy $\Ekininc=25$\,meV is marked by the dashed vertical line.}
\label{fig:rel_errors_scattering_states_1d}
\end{figure}

We repeat the same numerical experiment, but this time 
we compute scattering states for the whole range of energies 
$\Ekininc \in [10^{-3},10^3]$\,meV.
Simultaneously, we calculate (quasi)-exact reference solutions.
This is possible since the potential energy is a piecewise linear function 
\cite{LuFu89,VaGi96}.
The relative errors in the $\ell^2$-norm for different numerical methods
are depicted in Figure \ref{fig:rel_errors_scattering_states_1d} (left).
It turns out that $\mathrm{DTBC}_\mathrm{2nd}$ and $\mathrm{PML}_\mathrm{2nd}$ 
yield similar results for medium to large energies.
For very small energies, $\mathrm{DTBC}_\mathrm{2nd}$ performs significantly 
better than $\mathrm{PML}_\mathrm{2nd}$. Moreover, the
higher-order methods $\mathrm{PML}_\mathrm{4th}$ and $\mathrm{PML}_\mathrm{6th}$
yield much smaller errors than $\mathrm{DTBC}_\mathrm{2nd}$
for energies $\Ekininc \gtrsim 1$\,meV.
We remark that a smaller mesh size $\triangle x$ yields similar curves 
which are shifted downwards.

For zero potential energy $V=0$, the scattering state solutions to the 
stationary Schr\"o\-dinger equation are simple plane waves.
The relative errors, to which we will return later, are shown in 
Figure \ref{fig:rel_errors_scattering_states_1d} (right).

\begin{figure}[t]
  \centering
	\includegraphics[width=75mm]{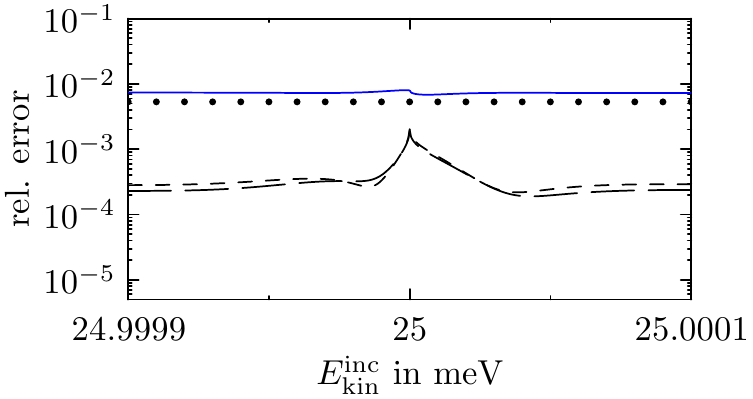}
  	\quad
  	\includegraphics[width=75mm]{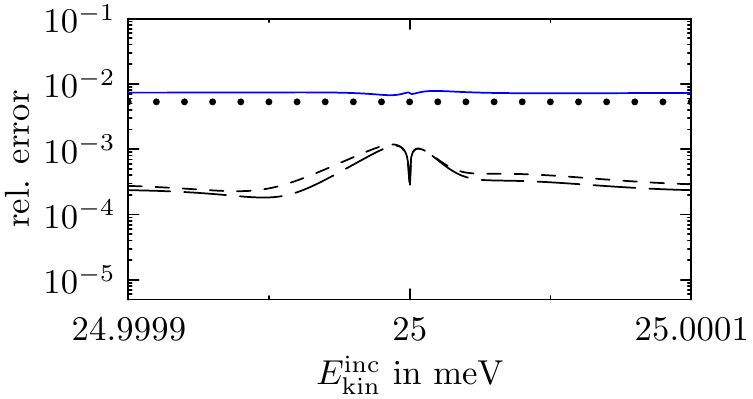}
  	\\
	\includegraphics[width=75mm]{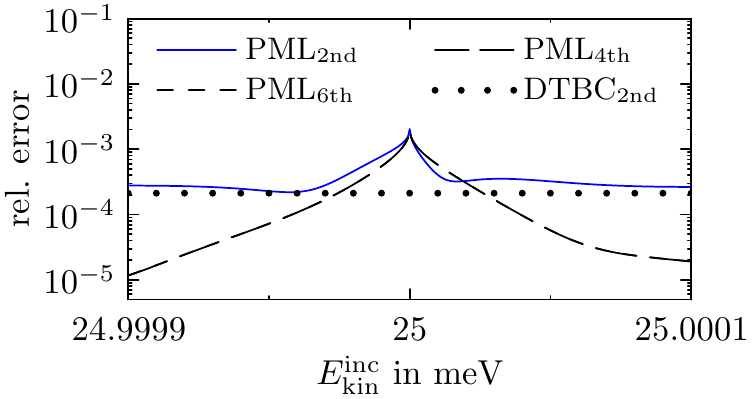}
  	\quad
	\includegraphics[width=75mm]{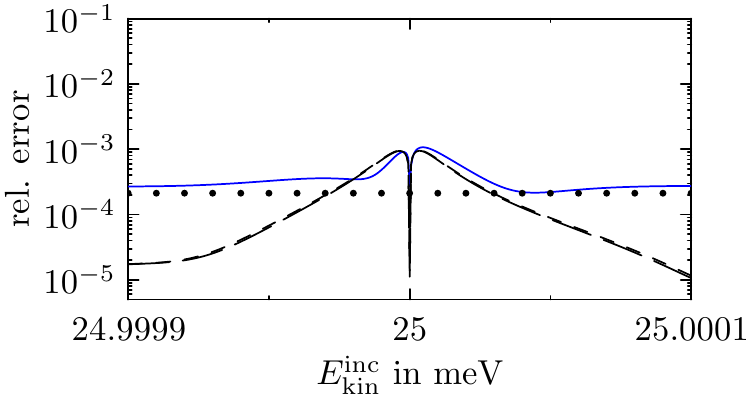}
  	\caption{Close-up view of the relative errors of the scattering states 
	as a function of the kinetic energy using Dirichlet
	(left column) or Neumann boundary conditions (right column) at the end points
	of the PML with mesh sizes $\triangle x=0.5$\,nm (top row) and 
	$\triangle x=0.1$\,nm (bottom row).}
  	\label{fig:rel_errors_pml_neumann_vs_dirichlet}
\end{figure}

Next, we compute scattering states for the ramp-like potential energy and 
$\Ekininc \in [24.9999,25.0001]$\,meV, 
i.e., for energies extremely close to the critical energy $\Ekininc=25$\,meV.
This time, we use Dirichlet or Neumann boundary conditions at the end points 
of the PML and employ two different mesh sizes $\triangle x=0.5$\,nm and 
$\triangle x=0.1$\,nm.  
The numerical errors corresponding to Dirichlet boundary conditions are 
depicted in the left column of Figure \ref{fig:rel_errors_pml_neumann_vs_dirichlet} 
and those corresponding to Neumann boundary conditions are shown in the right column.

For $\triangle x=0.5$\,nm (top row), we observe only a small perturbation 
of the numerical errors around the critical value. 
In case of $\mathrm{PML}_\mathrm{2nd}$, the effect is obscured completely 
by the error of the spatial discretization, which is present even if
transparent boundary conditions are used (compare $\mathrm{DTBC}_\mathrm{2nd}$).
For $\triangle x=0.1$\,nm (bottom row), the effect becomes more pronounced since
the spatial discretization error is greatly reduced.
However, the maximum error is essentially the same.
This holds also true for smaller $\triangle x$ and finer sampling of $\Ekininc$.
For $\Ekininc=25$\,meV the Neumann boundary condition is exact and hence, 
the numerical error is minimized. 

For that reason we prefer Neumann boundary 
conditions at the end points of the computational domain in all subsequent simulations.
Apart from this, both boundary conditions provide essentially the same results.
In summary, we found that the accuracy of the numerical methods using PML is reduced 
if the wave number of an incident wave function approaches zero.
However, in practice this effect is comparatively small.

\subsection{Wave packets}
\label{subsec:wave_packets_1d}

We consider the one-dimensional time-dependent Schr\"odinger equation,
\begin{equation}
  \label{eq:time_dependent_schroedinger_equation_1d}
  i \hbar \frac{\partial}{\partial t} \psi(x,t) = -\frac{\hbar^2}{2 m^\star} \frac{  
	\partial^2}{\partial x^2} \psi(x,t)
  + V(x,t) \psi, \quad \psi(\cdot,0) = \psi_0, \quad x \in \mathbb{R},\, t > 0.
\end{equation}
Under the assumptions that the initial wave function is compactly supported in 
$(0,L)$ and that the potential in the exterior domain vanishes, $V(x,t) = 0$ 
for $x \leq 0$ and $x \geq L$, $t \geq 0$, it is well known \cite{Arn01} 
that transparent boundary conditions at $x=0$ and $x=L$ read as
\begin{equation}
  \label{eq:tbc_zero_exterior_potential_1d}
  \frac{\pa \psi}{\pa x}(x,t) \big|_{x=0,L} 
	= \pm \sqrt{\frac{2 m^\star}{\pi \hbar}} e^{-i \pi/4}
  \frac{d}{dt}\int_0^t 
  \frac{\psi(x,\tau) \big|_{x=0,L}}{\sqrt{t-\tau}} \,d \tau.  
\end{equation}
Using \eqref{eq:tbc_zero_exterior_potential_1d} one can easily derive transparent 
boundary conditions for non-zero exterior potentials which are spatially constant 
but may change with time \cite{AnBe03}. As an example, we consider
$V(x,t) = V_r(t)$ for $x \geq L,\, t \geq 0$ and
$V_r(t) = -e U(t)$, where $U(t)$ is the applied voltage.
To get rid of the potential in the right lead, we define the gauge change
$$
  \widetilde{\psi}(x,t) := \exp\left(\frac{i}{\hbar}\int_0^t V_r(s)\,ds\right)\psi(x,t)
  \quad x \geq L,\; t \geq 0.
$$
This function solves the free Schr\"odinger equation and consequently, 
\eqref{eq:tbc_zero_exterior_potential_1d} yields a transparent boundary 
condition for a time-dependent exterior potential in the right lead: 
$$
  \frac{\partial \psi}{\partial x} (x,t)\big|_{x=L} 
	= -\sqrt{\frac{2 m^\star}{\pi \hbar}} e^{-i \pi / 4} 
  e^{-i \int_0	^t V_r(s)\,ds / \hbar}
  \frac{d}{dt} \int_0^t
  \frac{e^{i \int_0^t V_r(s)\,ds/\hbar}\psi(x,\tau)\big|_{x=L}}{\sqrt{t-\tau}} 
	\, d \tau.
$$

\medskip\noindent
{\bf Discrete transparent boundary conditions.}
The three-point finite-difference discretization \eqref{eq:D_2_2nd} applied to
the time-dependent Schr\"odinger equation 
\eqref{eq:time_dependent_schroedinger_equation_1d}
yields the semi-dis\-cre\-tized problem
\begin{equation}
  \label{eq:tdse_semi_discretized_problem_1d}
  \frac{d}{d t} \psi_j(t) = i \frac{\hbar}{2 m^\star} D_x^{2,\mathrm{2nd}} \psi_j(t) -  \frac{i}{\hbar} V(t) \psi_j(t)
  =:f(t, \psi_j(t)),
\end{equation}
which is solved by the Crank-Nicolson time-integration method
\begin{equation}
  \label{eq:Crank_Nicolson_time_integration}
  \psi^{(n+1)} = \psi^{(n)} + \triangle t \, 
  f\left( (n+1/2) \triangle t, \psi^{(n+1/2)} \right).
\end{equation}
Replacing $\psi_j^{(n+1/2)}$ by the average value $(\psi_j^{(n+1)}+\psi_j^{(n)})/2$
yields the well-known Crank-Nicolson scheme
\begin{equation}
  \label{eq:Crank_Nicolson_scheme_1d_2nd_order}
  \left(I- \frac{i \hbar \triangle t}{4 m^\star} D_x^2
  + \frac{i \triangle t}{2 \hbar} V_j^{(n+1/2)}\right) \psi_j^ {(n+1)}
  = \left(I + \frac{i \hbar \triangle t}{4 m^\star} D_x^2
  -\frac{i \triangle t}{2 \hbar} V_j^{(n+1/2)}\right) \psi_j^{(n)} 
\end{equation}
on the equidistant grid $x_j=j \triangle x$, $t_n=n\triangle t$
with $j \in \mathbb{Z}$ and $n \in \mathbb{N}_0$.
For zero exterior potentials, the corresponding DTBC at the left ($x_0=0$) 
and the right ($x_J=L$) contact are given as follows \cite{Arn98}:
\begin{subequations}
  \label{eq:dtbc_zero_exterior_potential_1d}
  \begin{align}
  \label{eq:hom_dtbc_1d_lhs}
  \psi_1^{(n+1)}-s^{(0)} \psi_0^{(n+1)} 
  &= \sum_{\ell=1}^{n} s^{(n+1-\ell)} \psi_0^{(\ell)} - \psi_1^{(n)}, \quad n\geq 0, 
  \\
  \label{eq:hom_dtbc_1d_rhs}
  \psi_{J-1}^{(n+1)} - s^{(0)} \psi_J^{(n+1)} 
  &= \sum_{\ell=1}^{n} s^{(n+1-\ell)} \psi_J^{(\ell)} 
	- \psi_{J-1}^{(n)}, \quad n \geq 0,
  \end{align}
\end{subequations}
with the convolution coefficients
$$
  s^{(n)}=\left( 1 - i \frac{R}{2} \right) \delta_{n,0}
  + \left( 1 + i \frac{R}{2} \right) \delta_{n,1} 
  + \alpha e^{-i n \varphi} \frac{P_n(\mu)-P_{n-2}(\mu)}{2n-1}
$$
and the abbreviations
$$
  R = \frac{4 m (\triangle x)^2}{\hbar \triangle t},\quad
  \varphi=\arctan\frac{4}{R}, \quad \mu=\frac{R}{\sqrt{R^2+16}}, 
  \quad \alpha = \frac{i}{2} \sqrt[4]{R^2(R^2+16)} e^{i \varphi / 2} .
$$
Here, $P_n$ denotes the $n$th-degree Legendre polynomial ($P_{-1}=P_{-2}=0$), 
and $\delta_{n,j}$ is the Kronecker symbol.
The Crank-Nicolson scheme \eqref{eq:Crank_Nicolson_scheme_1d_2nd_order} with 
the DTBC \eqref{eq:dtbc_zero_exterior_potential_1d}
yields an unconditionally stable discretization which is perfectly free of 
reflections \cite{Arn98, Arn01}.
The corresponding solution coincides exactly with the solution of the discrete 
whole space problem \eqref{eq:Crank_Nicolson_scheme_1d_2nd_order} restricted to
the grid $X_\mathrm{DTBC}$ defined in \eqref{eq:X_dtbc}.
Time-dependent exterior potentials may be included like in the continuous case 
described above (see \cite{MJK12} for details).

\medskip\noindent
{\bf Perfectly Matched Layers.}
We apply the coordinate transformation \eqref{eq:complex_coordinate_transformation} 
also in the transient case, which yields
the time-dependent Sch\"odinger-PML equation
\begin{equation}
  \label{eq:time_dependent_schroedinger_pml_equation_1d}
  i \hbar \frac{\partial}{\partial t} \psi(x,t) = 
  -\frac{\hbar^2}{2 m^\star} c(x) \frac{\partial}{\partial x} 
	\left( c(x) \frac{\partial}{\partial x} \psi(x,t) \right) + V(x,t) \psi(x,t),
\end{equation}
where $c(x)$ is defined in \eqref{eq:pml_complex_function_c}.
Accordingly, the semi-discretized problem
is given by \eqref{eq:tdse_semi_discretized_problem_1d},
wherein $D_x^{2,\mathrm{2nd}}$ is replaced by
$\tilde{D}_x^2 \in \{\tilde{D}_x^{2,\mathrm{2nd}}, \tilde{D}_x^{2,\mathrm{4th}}, 
\tilde{D}_x^{2,\mathrm{6th}}\}$ (see \eqref{eq:D_tilde_1d}).
The spatial grid $X_\mathrm{PML}$ is the same 
as in the stationary case (see \eqref{eq:X_pml}).
Using the Crank-Nicolson time-integration method gives 
\eqref{eq:Crank_Nicolson_scheme_1d_2nd_order} but with the modified 
spatial differential operator $\tilde{D}_x^2$.

Alternatively, we solve this problem
via the classical (explicit) Runge Kutta method.
In quantum mechanics simulations, this method is used very rarely.
This is probably because the resulting spatio-temporal discretization 
is only conditionally stable. Moreover, the norm of the wave function
is not a conserved quantity, i.e., in general 
$\|\psi^{(n+1)}\|_{\ell^2}^2 = \|\psi^{(n)}\|_{\ell^2}^2$ 
does not hold exactly.
We address these issues in a simple numerical experiment where we solve the 
ordinary time-dependent Schr\"odinger equation (without PML) for 
a harmonic oscillator potential $V(x) = m^\star \omega_\ast x^2 / 2$ with
$\omega_\ast = 0.25 \times 10^{14}\mathrm{\,s}^{-1}$. 
As a reference solution, we consider a so called coherent state \cite{Ga09},
$$
  \psi(x,t) = \left( \frac{m^\star \omega_\ast}{\pi \hbar} \right)^{1/4}
  \exp\left(-\frac{m^\star \omega_\ast}{2 \hbar}
	\Big(x^2-2 x x_0 e^{-i \omega_\ast t}+\frac{x_0^2}{2} e^{-2 i \omega_\ast t}
	+\frac{x_0^2}{2}\Big)-\frac{i}{2} \omega_\ast t\right),
$$
where $x_0=10$\,nm denotes the expectation value of the particle at $t=0$.
The computational domain extends from $-50$\,nm to $50$\,nm.
Near the boundaries, the wave function is zero (to numerical precision)
and hence it is reasonable to employ homogeneous Dirichlet boundary conditions.
The initial state $\psi(x,0)$ is propagated for $100\,000$ time steps using the
Runge-Kutta scheme.
The spatial derivative is approximated by $D_x^{2,\mathrm{2nd}}$, 
$D_x^{2,\mathrm{4th}}$, or $D_x^{2,\mathrm{6th}}$.
Using $\triangle x=0.5$\,nm and $\triangle t=0.1$\,fs, we obtain
the deviation from the initial mass, 
$\big| \| \psi^{(100\,000)} \|_{\ell^2} \big/ \| \psi^{(0)} \|_{\ell^2} - 1 \big| 
\lesssim 5.9 \times 10^{-11}$,
independent of the spatial discretization.
Moreover, the relative errors of the final wave functions are given by
$2.19 \times 10^{-1}$, $6.56 \times 10^{-4}$, and $6.42 \times 10^{-6}$ 
corresponding to $D_x^{2,\mathrm{2nd}}$, $D_x^{2,\mathrm{4th}}$, and 
$D_x^{2,\mathrm{6th}}$, respectively.
In the recent paper \cite{CaCa13}, the stability bounds 
$\triangle t < \kappa (\triangle x)^2$
for the linear and nonlinear Schr\"odinger equation were derived 
for second- and fourth-order spatial discretizations.
Remarkably, the constants $\kappa$ stated in \cite{CaCa13} agree 
with our experimental findings, even if we solve the Schr\"odinger-PML equation.
In case of the sixth-order spatial discretization, our numerical experiments suggest 
that $\kappa$ needs to be adapted slightly by a factor $3/4$ compared to $\kappa$ 
in the fourth-order discretization. In fact, 
$\kappa = 9 m^\star/(8 \sqrt{2} N \hbar)$,
where $N$ denotes the space dimension,
gives almost sharp bounds for the simulations presented in this paper.

\medskip\noindent
{\bf Simulations.}
We solve the time-dependent Schr\"odinger equation
\eqref{eq:time_dependent_schroedinger_equation_1d} for zero potential energy.
As a reference solution, we choose a superposition of three Gaussian wave packets,
\begin{equation}
  \label{eq:gaussian_1d}
  \xi_p(x,t) = \left( 1 + i \frac{t}{\tau} \right)^{-1/2} 
	\exp\bigg[\left( 1 + i \frac{t}{\tau} \right)^{-1}
	\bigg(-\left( \frac{x-x_0}{2 \sigma} \right)^2 + i k_p (x-x_0) -i \sigma^2 k_p^2 
	\frac{t}{\tau}\bigg)\bigg],
\end{equation}
where $\tau=2 m^\star \sigma^2/\hbar$ and $k_p=\sqrt{2m^\star E_p}/\hbar$, $p=1,2,3$.
At $t=0$, each Gaussian is centered around $x_0=L/2$.
Using $\sigma=7.5$\,nm, the initial wave packet is practically zero outside 
the device domain $[0,L]=[0,120]$\,nm.
The average energy of the first, second, and third Gaussian is $E_1=0$\,meV, 
$E_2=25$\,meV, and $E_3=75$\,meV, respectively.
Therefore, the reference solution is a superposition of propagating plane waves 
$\exp(i k x -i \omega t)$ of very low to high energies 
$E=(\hbar k)^2 / (2 m^\star)$ with the wave frequency $\omega=E/\hbar$.
We use the same spatial mesh size $\triangle x=0.5$\,nm 
and the same parameters of the PML as for the scattering state simulations. 
The time step size is given by $\triangle t=0.1$\,fs.

\begin{figure}[thb]
	\centering
  	\includegraphics[width=150mm]{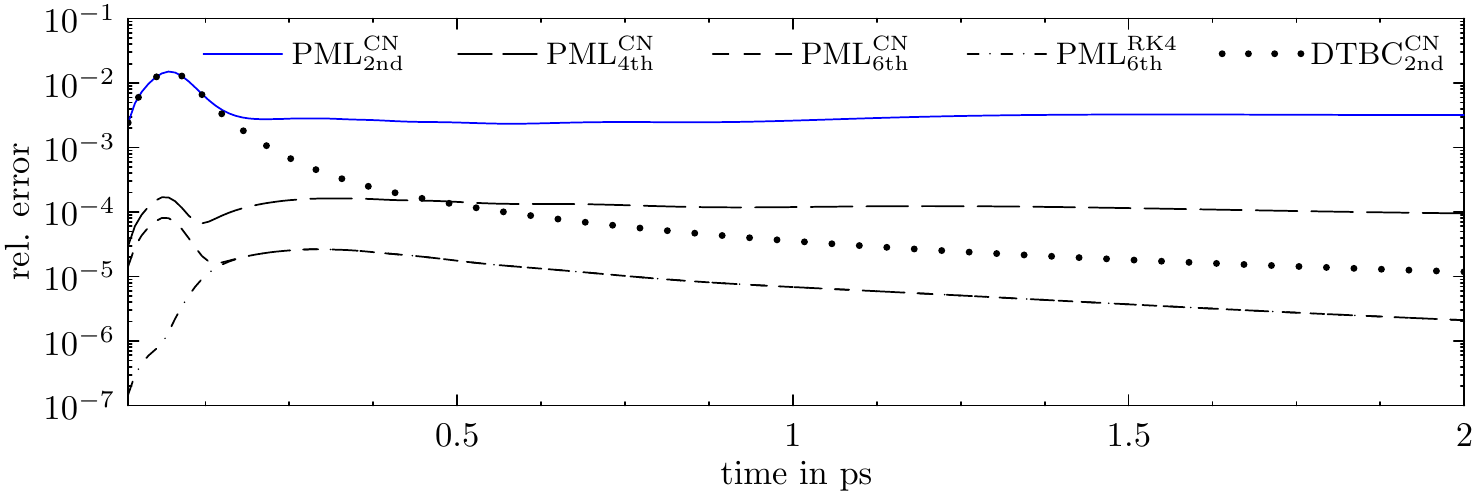}
  	\caption{Relative errors versus time for different numerical methods 
	corresponding to the simulation using three Gaussian wave packets.}
  \label{fig:rel_errors_wave_packets_1d}
\end{figure}

Figure \ref{fig:rel_errors_wave_packets_1d} shows the relative $\ell_2$ errors
of simulations using DTBC or PML.
At the beginning of the simulation, the two second-order methods yield 
similar results. During this phase, the numerical error is dominated by the 
fast traveling parts of the wave packet, since their oscillations in space and 
time are difficult to handle by low order methods.
Shortly afterwards, the relative error of the PML solution is stabilizing 
around the value $3\times 10^{-3}$.
This agrees with the numerical errors presented in 
Figure \ref{fig:rel_errors_scattering_states_1d} (right) since at later times, 
the wave packet can be thought of as a superposition of primarily low energy 
plane wave scattering states.
In contrast, the numerical error of the DTBC solution decreases continuously.

Figure \ref{fig:rel_errors_wave_packets_1d} also shows the numerical errors 
according to higher-order methods.
For example, a sixth-order spatial discretization in combination with PML 
reduces the maximum relative error by more than two orders of magnitude.
The maximum error can be reduced further by replacing the second-order 
Crank-Nicolson time-integration method with the fourth-order Runge-Kutta method.
In this case, even the fast temporal oscillations in the beginning of the simulation
are resolved with high accuracy.


\subsection{Time-dependent incoming waves}
\label{subsec:time_dependent_incoming_waves_1d}

Before we turn our attention to transient scattering state simulations, we 
explain how to realize an incoming wave at a device boundary.
As an example, we consider the free time-dependent Schr\"odinger equation, 
where an incident plane wave 
$$
  \psi^{\mathrm{inc}}(x,t) = \exp(i k x - i \omega t), \quad x \leq 0,\quad  
	\omega = \hbar k^2 / (2 m^\star),
$$
is prescribed at $x=0$\,nm.
The transparent boundary condition at $x=L$ is as in 
\eqref{eq:tbc_zero_exterior_potential_1d}.
However, at $x=0$ we need to prescribe an inhomogenous transparent boundary condition 
which follows if we apply \eqref{eq:tbc_zero_exterior_potential_1d} to the wave 
$\psi-\psi^\mathrm{inc}$  \cite{AnBe03,Arn01}:
\begin{equation}
  \label{eq:inhomogeneous_tbc_zero_exterior_potential_1d}
  \frac{\partial}{\partial x} \left( \psi(x,t) - \psi^\mathrm{inc}(x,t) \right) 
	\big|_{x=0} = \sqrt{\frac{2 m^\star}{\pi \hbar}} e^{-i \pi/4}
  \frac{d}{dt}\int_0^t 
  \frac{ \psi(0,\tau) - \psi^\mathrm{inc}(0,\tau) }{\sqrt{t-\tau}} \,d \tau.  
\end{equation}
To avoid a discontinuity at $(x,t)=(0,0)$, we assume compatibility of the initial
and boundary data at this point (see \cite{Arn01} for details).

\medskip\noindent
{\bf Discrete transparent boundary conditions.}
The discrete analogue of \eqref{eq:inhomogeneous_tbc_zero_exterior_potential_1d}
for the Crank-Nicolson scheme \eqref{eq:Crank_Nicolson_scheme_1d_2nd_order}
follows by replacing $\psi_j^{(n)}$ by $\psi_j^{(n)}-\phi_j^{(n)}$
in \eqref{eq:hom_dtbc_1d_lhs},
$$
  \psi_1^{(n+1)}-s^{(0)} \psi_0^{(n+1)} 
  = \sum_{\ell=1}^{n} s^{(n+1-\ell)} \left( \psi_0^{(\ell)} - \phi_0^{(\ell)}\right)
  - \left( \psi_1^{(n)} - \phi_1^{(n)}\right)
  + \phi_1^{(n+1)} - s^{(0)} \phi_0^{(n+1)}.
$$
Here,
\begin{equation}
  \label{eq:time_dependent_discrete_plane_wave_1d}
  \phi_j^{(n)} = \exp(i k x_j - i \omega n \triangle t)
\end{equation}
represents an incoming discrete plane wave, i.e.,
$k$ is related to $E=\Ekininc$ according to 
the discrete $E$--$k$--relation \eqref{eq:discrete_E_k_relation},
and the wave frequency is given by the discrete $E$--$\omega$--relation
\begin{equation}
  \label{eq:discrete_E_omega_relation}
  \omega = \frac{2}{\triangle t} \arctan\left(\frac{E \triangle t}{2 \hbar}\right),
\end{equation}
which is the discrete analogue of $\omega = E / \hbar$ (see \cite{Arn01}).

\medskip\noindent
{\bf Perfectly Matched Layers.}
We solve the transient Schr\"odinger-PML equation 
\eqref{eq:time_dependent_schroedinger_pml_equation_1d}, where an incoming 
time-dependent plane wave is prescribed at one of the device boundaries.
The incoming wave is realized analogously to the case of the 
stationary scattering state simulations.
As an example, we consider the Crank-Nicolson scheme
\eqref{eq:Crank_Nicolson_scheme_1d_2nd_order} with $D_x^2$ replaced by
$\tilde{D}_x^2 = \tilde{D}_x^{2, \mathrm{2nd}}$.
In the vicinity of the device boundary $x_{j_0}=0$,
the potential energy is zero, but more importantly, the PML is inactive.
The wave function in the left contact is a superposition of an incoming 
and a reflected wave.
However, the incoming wave \eqref{eq:time_dependent_discrete_plane_wave_1d}
is eliminated in the left lead and hence $\psi_j^{(n)}$ represents the 
reflected wave for $j<j_0$. 
Accordingly, the finite-difference equations for $j_0-1$ and $j_0$ need 
to be modified as follows:
\begin{align}
  \label{eq:modifications_inc_wave_pml_1d}
  \begin{split}
  \psi_{j_0-1}^{(n+1)} 
  &- \frac{i \hbar \triangle t}{4 m^\star (\triangle x)^2}
  \Big[\psi_{j_0-2}^{(n+1)} - 2 \psi_{j_0-1}^{(n+1)} + \big( \psi_{j_0}^{(n+1)} 
	- \phi_{j_0}^{(n+1)}\big) \Big] \\
  &= \psi_{j_0-1}^{(n)} + \frac{i \hbar \triangle t}{4 m^\star (\triangle x)^2}
  \Big[ \psi_{j_0-2}^{(n)} - 2 \psi_{j_0-1}^{(n)} + \big( \psi_{j_0}^{(n)} 
	- \phi_{j_0}^{(n)} \big) \Big], \\
  \psi_{j_0}^{(n+1)} 
  &-\frac{i \hbar \triangle t}{4 m^\star (\triangle x)^2}
  \Big[\big( \psi_{j_0-1}^{(n+1)} + \phi_{j_0-1}^{(n+1)} \big) - 2 \psi_{j_0}^{(n+1)} 
	+ \psi_{j_0+1}^{(n+1)}\Big] \\
  &= \psi_{j_0}^{(n)} + \frac{i \hbar  \triangle t}{4 m^\star (\triangle x)^2}
  \Big[ \big( \psi_{j_0-1}^{(n)} + \phi_{j_0-1}^{(n)} \big) - 2 \psi_{j_0}^{(n)} 
	+ \psi_{j_0+1}^{(n)} \Big].
  \end{split}
\end{align}
Since $\phi_{j_0}^{(n)}$ and $\phi_{j_0-1}^{(n)}$ are known for all 
$n \in \mathbb{N}_0$, we collect these values,
$$
  b_j^{(n)} =
  \begin{cases}
  \begin{aligned}
  -\frac{i \hbar \triangle t}{4 m^\star (\triangle x)^2)}
  \big(\phi_{j_0}^{(n+1)} + \phi_{j_0}^{(n)}\big) 
	&\quad \text{for } j = j_0-1, \\
  +\frac{i \hbar \triangle t}{4 m^\star (\triangle x)^2)}
  \big(\phi_{j_0-1}^{(n+1)} + \phi_{j_0-1}^{(n)}\big) 
	&\quad \text{for } j = j_0, \\
  0 &\quad \textrm{else}
  \end{aligned}
  \end{cases}
$$
on the right-hand side of \eqref{eq:Crank_Nicolson_scheme_1d_2nd_order}.
In case $\tilde{D}_x^2=\tilde{D}_x^{2,\mathrm{4th}}$ and 
$\tilde{D}_x^2=\tilde{D}_x^{2,\mathrm{6th}}$, we proceed in a similar way.
However, due to the extended finite-difference stencils, four or six 
finite-difference equations need to be modified accordingly.
If the semi-discretized problem is solved via the Runge-Kutta method,
the incoming wave needs to be taken into account 
at each of the four intermediate Runge-Kutta time-steps. 

\begin{figure}
	\centering
	\includegraphics[width=75mm]{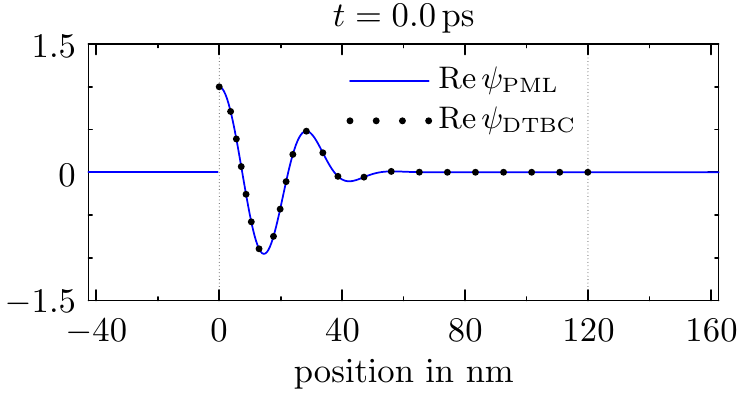}
  	\hspace{2.5mm}
	\includegraphics[width=75mm]{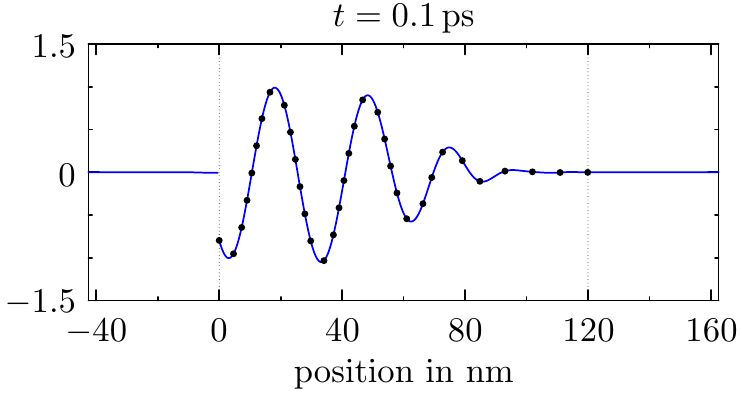}
  	\caption{An incoming plane wave with energy $\Ekininc = 25$\,meV 
	at the left contact using PML (solid line) and DTBC (dotted line).}
\label{fig:incoming_plane_wave_1d}
\end{figure}

\medskip\noindent
{\bf Simulations.}
An incoming plane wave is depicted in Figure~\ref{fig:incoming_plane_wave_1d}
at $t=0$\,ps and $t=0.1$\,ps. 
The potential energy is zero everywhere and
the kinetic energy of the incoming electrons amounts to $\Ekininc=25$\,meV.
We note that it would take quite a long time (compared to the underlying time scale)
before the wave function becomes approximately stationary.

\subsection{Transient scattering states}
\label{subsec:transient_scattering_1d}

A transient scattering state simulation describes a quantum device in which 
a continuously incoming plane wave is prescribed at one of the device contacts.
Moreover, the potential energy is allowed to change with time.
The lead potentials need to be spatially constant but may depend on time too.
In contrast to the simulation depicted in Figure \ref{fig:incoming_plane_wave_1d}, 
one is typically not interested in the initial transient phase.
Instead of waiting for the wave function to become stationary for the first time, 
it is preferable to initialize the simulation with a scattering state.
This situation is considered in \cite{Arn01} where the potential energy is 
switched instantaneously. An extension to continuously time-variable potentials 
can be found in \cite{MJK12}.

As an example, we consider the time-dependent Sch\"odinger equation with the 
ramp-like potential \eqref{eq:ramp_like_potential}, 
where the applied voltage $U$ is assumed to be time-dependent with $U(t) = U_0$
for $t\le 0$. The initial wave function 
is given by the scattering state solution $\phi$ of the stationary Schr\"odinger 
equation for the potential energy according to $U_0$.
As before, we consider electrons injected at the left contact with energy 
$E=\Ekininc$. For the sake of completeness, we state the corresponding 
boundary conditions \cite{Arn01,MJK12}:
\begin{subequations}
  \label{eq:inhom_tbc_time_dependend_ext_potential_1d}
  \begin{align}
  \label{eq:inhom_tbc_time_dependend_ext_potential_1d_lhs}
  \begin{split}
  &\frac{\partial}{\partial x} \left[ \psi(x,t)-e^{-i E t / \hbar}
	\phi(x)\right]\Big|_{x=0} \\
  &\quad\quad = +\sqrt{\frac{2 m^\star}{\pi \hbar}} e^{-i \pi / 4} 
  \frac{d}{dt} \int_0^t
  \frac{\psi(0,\tau) - e^{-i E \tau / \hbar} \phi(0)}{\sqrt{t-\tau}} \, d \tau, 
  \end{split} \\
  \label{eq:inhom_tbc_time_dependend_ext_potential_1d_rhs}
  \begin{split}
  &\frac{\partial}{\partial x} 
  \left[ e^{i \int_0	^t V_r(s)\,ds / \hbar} \psi(x,t) -e^{-i (E-V_r(0)) t / \hbar} 
	\phi(x)\right]\Big|_{x=L} \\
  &\quad\quad =
  -\sqrt{\frac{2 m^\star}{\pi \hbar}} e^{-i \pi / 4} 
  \frac{d}{dt} \int_0^t
  \frac{e^{i \int_0^\tau V_r(s) \, ds / \hbar} \, \psi(L,\tau) - e^{-i(E-V_r(0))
	\tau/\hbar} \phi(L) }{\sqrt{t-\tau}} \, d \tau.
  \end{split}
  \end{align}
\end{subequations}
Here, $V_r(t) = -e U(t)$ denotes the potential energy in the right lead.

\medskip\noindent
{\bf Discrete transparent boundary conditions.}
The discrete analogue of \eqref{eq:inhom_tbc_time_dependend_ext_potential_1d_lhs}
follows by replacing $\psi_j^{(n)}$ by $\psi_j^{(n)}-\beta^{(n)} \phi_j$ 
in \eqref{eq:hom_dtbc_1d_lhs},
\begin{equation}
  \label{eq:inhom_dtbc_1d_lhs}
  \begin{aligned}
  (\psi_1^{(n+1)} -& \beta^{(n+1)} \phi_1)
  -s^{(0)} (\psi_0^{(n+1)} - \beta^{(n+1)} \phi_0) \\
  &= \sum_{\ell=1}^{n} s^{(n+1-\ell)} (\psi_0^{(\ell)} - \beta^{(\ell)} \phi_0)
  - ( \psi_1^{(n)} - \beta^{(n)} \phi_1 ),
  \end{aligned}
\end{equation}
where $\phi_j$ is a solution of the discrete scattering state problem outlined above.
The discretization of the gauge-change term,
\begin{equation*}
  \beta^{(n)}=\exp \left( 
  - 2 i n \arctan \left( \triangle t E / (2 \hbar) \right )\right)
  \approx \exp(-i E t / \hbar),
\end{equation*}
is consistent with the underlying Crank-Nicolson time-integration method 
\cite{MJK12}. Similarly, the discrete analogue of 
\eqref{eq:inhom_tbc_time_dependend_ext_potential_1d_rhs} follows by replacing
$\psi_j^{(n)}$ by $\epsilon^{(n)} \psi_j^{(n)} - \gamma^{(n)} \phi_j$
in \eqref{eq:hom_dtbc_1d_rhs}:
\begin{equation}
  \begin{aligned}
  \label{eq:inhom_dtbc_1d_rhs}
  (\epsilon^{(n+1)} \psi_{J-1}^{(n+1)} &- \gamma^{(n+1)} \phi_{J-1}) 
  - s^{(0)} (\epsilon^{(n+1)} \psi_J^{(n+1)} - \gamma^{(n+1)} \phi_J) \\
  &= \sum_{\ell=1}^{n} s^{(n+1-\ell)} (\epsilon^{(\ell)} \psi_J^{(\ell)} 
	- \gamma^{(\ell)} \phi_J) 
  - (\epsilon^{(n)} \psi_{J-1}^{(n)} - \gamma^{(n)} \phi_{J-1}).
  \end{aligned}
\end{equation}
The gauge-change terms are approximated via
\begin{equation}
  \label{eq:gauge_change_gamma_and_epsilon}
  \begin{aligned}
  \gamma^{(n)} &= \exp (2 i n ( \arctan (\triangle t V_r^{(0)} / (2 \hbar) ) 
  - \arctan ( \triangle t E / ( 2 \hbar ) ) ))\\
  &\approx \exp(-i (E-V_r(0)) t / \hbar),\\
	\epsilon^{(n)} &= \exp \bigg(\frac{i}{\hbar} \sum_{\ell=0}^{n-1} 
  \arctan ( \triangle t V_r^{(\ell+1/2)} ) \bigg)
  \approx \exp\bigg(\frac{i}{\hbar} \int_0^t V_r(s)\,ds\bigg),
  \end{aligned}
\end{equation}
which is also compatible with the underlying Crank-Nicolson method.
We note that the Crank-Nicolson scheme along with the inhomogeneous DTBC
\eqref{eq:inhom_dtbc_1d_lhs} and \eqref{eq:inhom_dtbc_1d_rhs}
is still perfectly free of spurious reflections.
In fact, \eqref{eq:inhom_dtbc_1d_lhs} and \eqref{eq:inhom_dtbc_1d_rhs} yield an
exact truncation of the discrete whole-space problem.

\medskip\noindent
{\bf Perfectly Matched Layers.}
No further steps are necessary to realize a transient scattering state 
simulation using PML. Given the potential energy at $t=0$ and the kinetic energy 
$\Ekininc$ of the electrons injected at the device contact,
we compute a scattering state solution of the discrete 
stationary Sch\"odinger-PML equation.
This scattering state serves as the initial state of the transient problem where a
time-dependent incoming plane wave \eqref{eq:time_dependent_discrete_plane_wave_1d} 
is prescribed at the device contact.

\begin{figure}
  \centering
  \includegraphics[width=75mm]{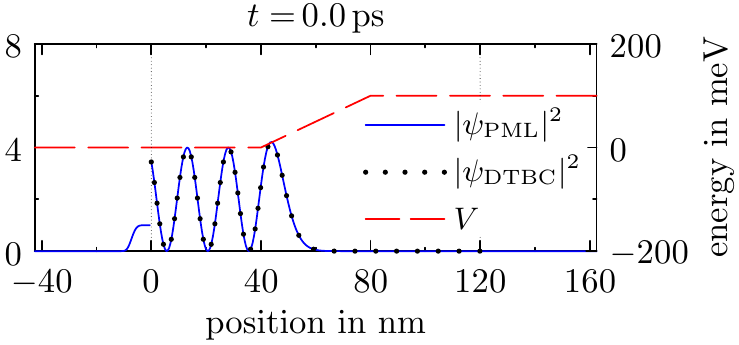}
  \includegraphics[width=75mm]{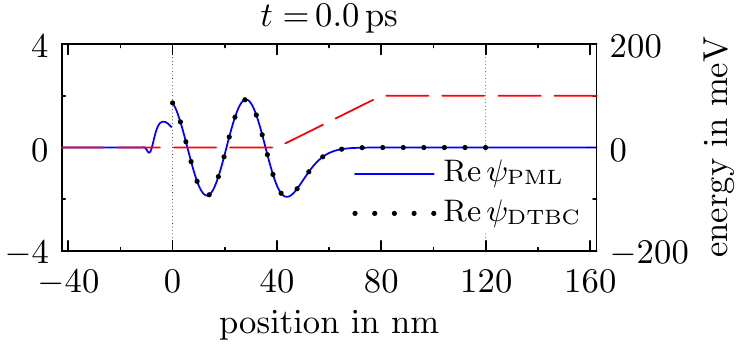} \\
  \vspace{2.0mm}
  \includegraphics[width=75mm]{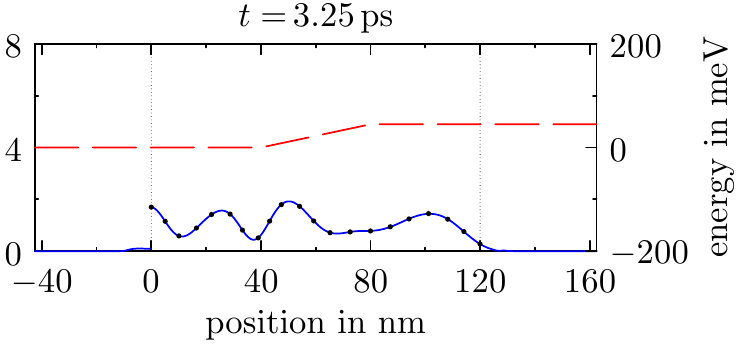}
  \includegraphics[width=75mm]{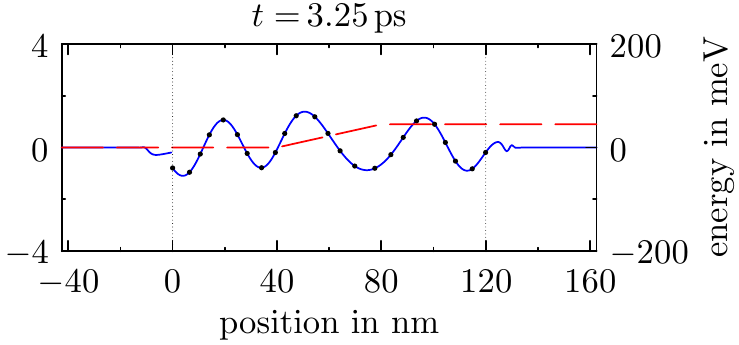} \\
  \vspace{2.0mm}
  \includegraphics[width=75mm]{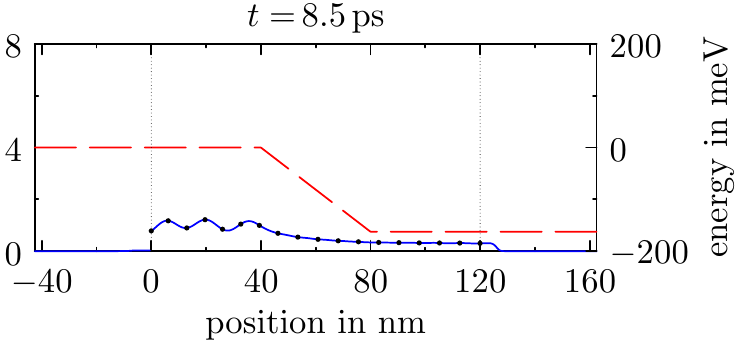}
  \includegraphics[width=75mm]{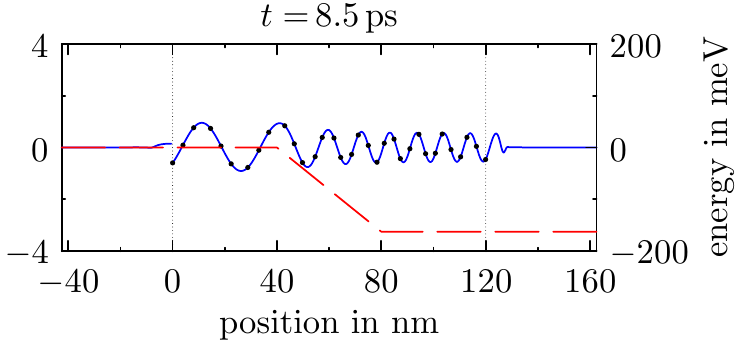} \\
  \vspace{2.0mm}
  \includegraphics[width=75mm]{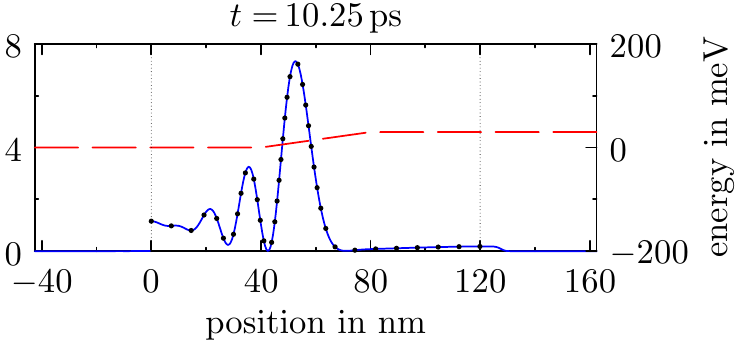}
  \includegraphics[width=75mm]{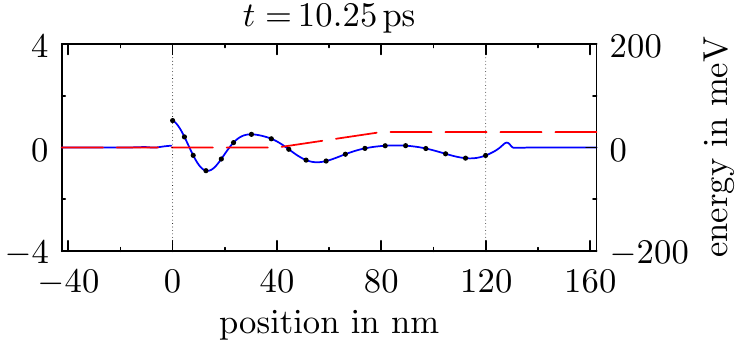} \\
  \vspace{2.0mm}
  \includegraphics[width=75mm]{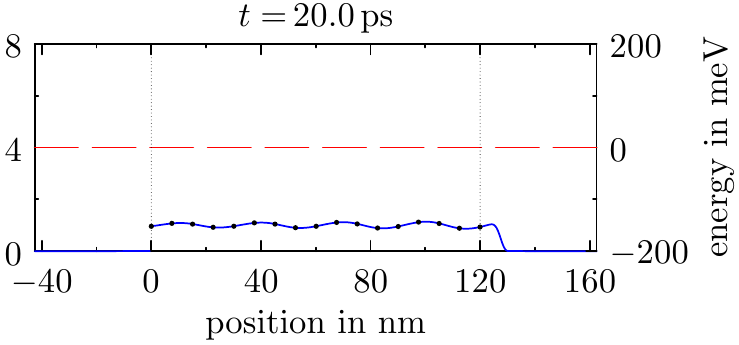}
  \includegraphics[width=75mm]{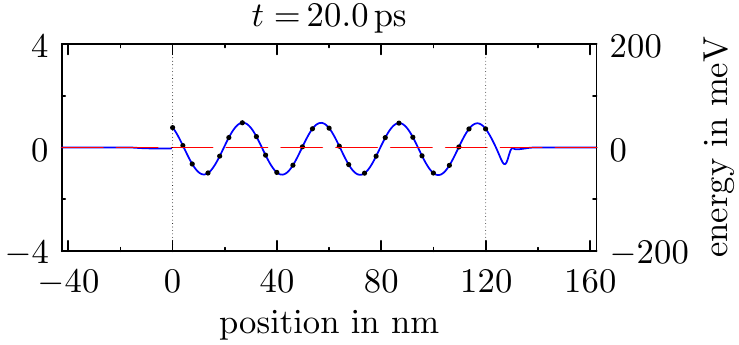}
  \caption{
  Transient scattering state $|\psi|^2$ (left column) and $\mbox{Re}\,\psi$ (right
	column) at selected times for a continuously incoming   
	plane wave prescribed at the left contact.
  The incoming plane wave corresponds to electrons with a kinetic energy of 
	$\Ekininc = 25$\,meV traveling to the right.}
  \label{fig:transient_scattering_1d_selected_times}
\end{figure}

\begin{figure}
  	\centering
  	\includegraphics[width=75mm]{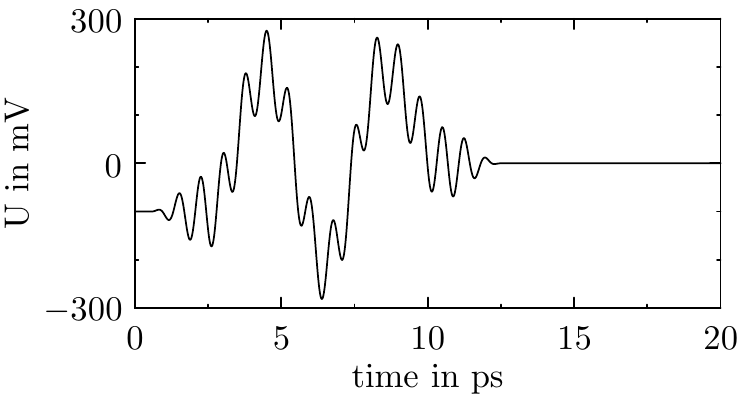}
  	\includegraphics[width=75mm]{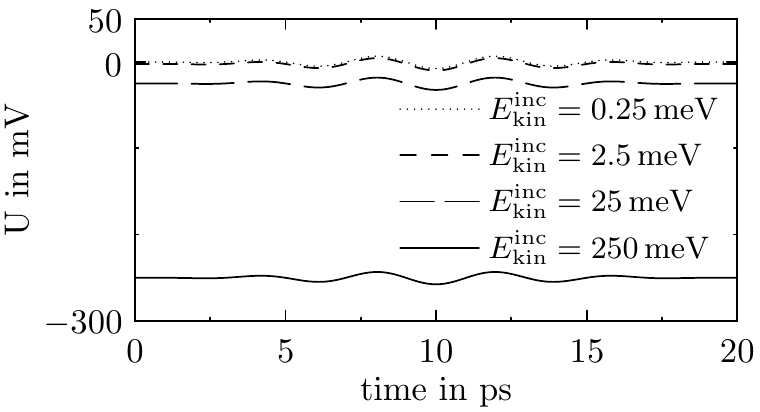} \\
	\includegraphics[width=75mm]{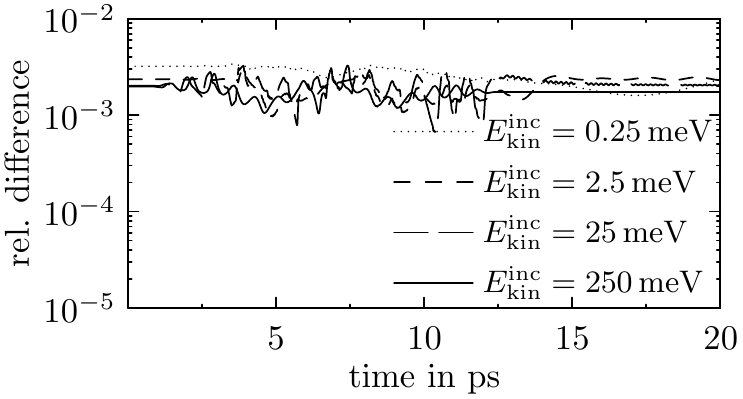}
  	\includegraphics[width=75mm]{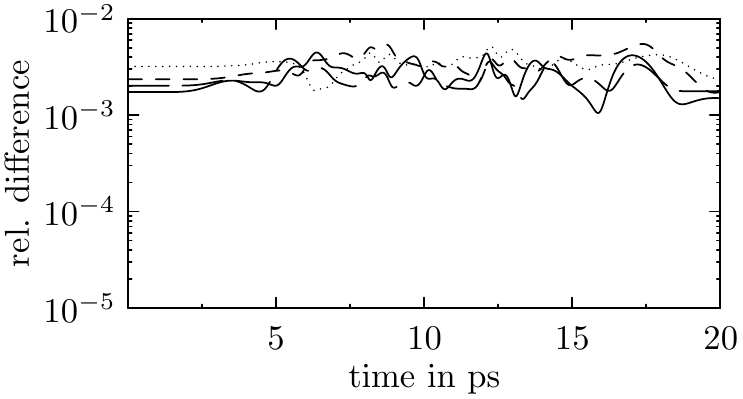} \\
	\includegraphics[width=75mm]{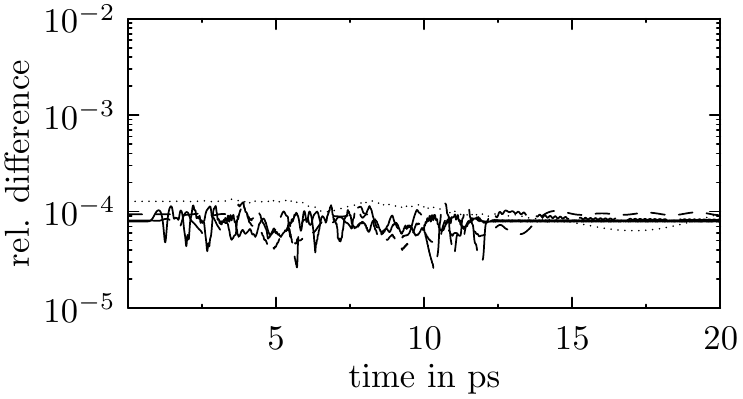}
	\includegraphics[width=75mm]{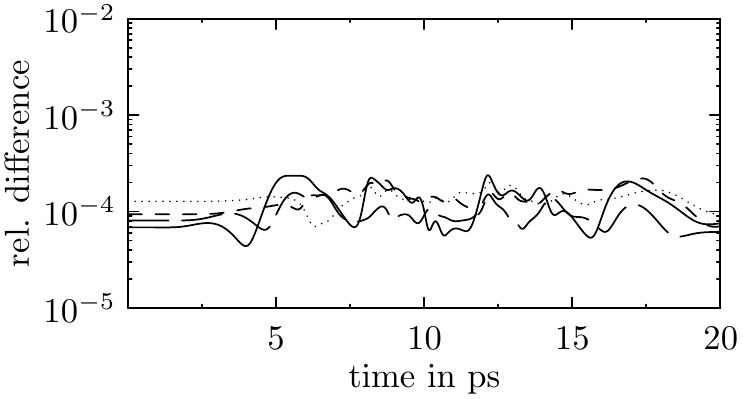}
	\caption{
	{\em Left column:} 
	Relative differences between the solutions computed with DTBC 
  	and PML for $\triangle x=0.5$\,nm (center) or
	$\triangle x=0.1$\,nm (bottom) using the applied voltage $U(t)$ (top).
  	The experiment is carried out four times corresponding to different values of
  	$\Ekininc=0.25$\,meV, $2.5$\,meV, $25$\,meV or $250$\,meV.
	{\em Right column:} 
	Relative differences between the solutions computed with DTBC 
	and PML for $\triangle x=0.5$\,nm (center) or
	$\triangle x=0.1$\,nm (bottom) using different applied voltages (top).
	The experiment is carried out four times corresponding to different values of
  	$\Ekininc=0.25$\,meV, $2.5$\,meV, $25$\,meV or $250$\,meV.
  	In each simulation the applied voltage oscillates slowly around a different
  	critical value
  	$U=-0.25$\,mV, $-2.5$\,mV, $-25$\,mV or $-250$\,mV.
}
	\label{fig:U_of_time_and_rel_differences_dtbc_pml_2nd_cn}
\end{figure}

\medskip\noindent
{\bf Simulations.}  
Figure \ref{fig:transient_scattering_1d_selected_times} shows the time evolution 
of a transient scattering state.
The incoming plane wave at the left contact represents electrons with a kinetic 
energy of $\Ekininc=25$\,meV traveling to the right.
The applied voltage $U(t)$ is illustrated in Figure 
\ref{fig:U_of_time_and_rel_differences_dtbc_pml_2nd_cn} (top left).
Accordingly, we initialize the simulation with the scattering state 
corresponding to an applied voltage of $U_0 = -100$\,mV.
Until $t=0.5$\,ps we keep the applied voltage constant and hence
the numerical solution remains the same.
More precisely, $|\psi|^2$ remains the same while $\Re \psi$ oscillates with time.

From $t=0.5$\,ps up to $t=12.5$\,ps, the applied voltage varies corresponding 
to a medium oscillation with a large amplitude and a fast oscillation 
with a smaller amplitude. 
As a result, the wave function shows a wild behavior as indicated in the 
second, third, and fourth row of Figure 
\ref{fig:transient_scattering_1d_selected_times}.
For times $t\ge 12.5$\,ps, the voltage is kept constant at $U=0$\,mV.
The last row of Figure \ref{fig:transient_scattering_1d_selected_times} shows that
even at $t=20$\,ps, the wave function has not become perfectly stationary again.
The time evolution of the scattering state is available as a movie at
\url{http://www.asc.tuwien.ac.at/~juengel}.

The simulation described above was carried out using the second-order
Crank-Nicolson scheme with DTBC or PML with the time step size of 
$\triangle t=0.1$\,fs.
Inside the device domain both discretizations coincide exactly.
Since DTBC represent an exact truncation of the discrete whole-space problem, 
the relative difference of both methods results from the PML.
We repeat the experiment four times corresponding to different values of
the kinetic energy $\Ekininc=0.25$\,meV, $2.5$\,meV, $25$\,meV or $250$\,meV.
The relative differences as a function of time are depicted 
in the left column of Figure \ref{fig:U_of_time_and_rel_differences_dtbc_pml_2nd_cn} 
for $\triangle x = 0.5$\,nm (center) or $\triangle x=0.1$\,nm (bottom).

In another four simulations, we let the applied voltage oscillate slowly
around the critical values $U=-0.25$\,mV, $-2,5$\,mV, $-25$\,mV or $-250$\,mV
as depicted in the right column of Figure 
\ref{fig:U_of_time_and_rel_differences_dtbc_pml_2nd_cn} (top).
In this way, we trigger waves of arbitrary low energy in the right lead.
The relative differences are depicted in the right column of
Figure \ref{fig:U_of_time_and_rel_differences_dtbc_pml_2nd_cn} 
for $\triangle x = 0.5$\,nm (center) or $\triangle x=0.1$\,nm (bottom).
It can be seen that the relative differences are alwas localized around a value of
$3\times 10^{-3}$, showing that PML can handle even the extreme cases.




\section{Two-dimensional quantum waveguide simulations}
\label{sec:quantum_waveguide_simulations}

\subsection{Scattering states in quantum waveguides}
\label{subsec:scattering_states_waveguide}

\begin{figure}[ht]
  	\begin{centering}
  	\includegraphics[width=150mm]{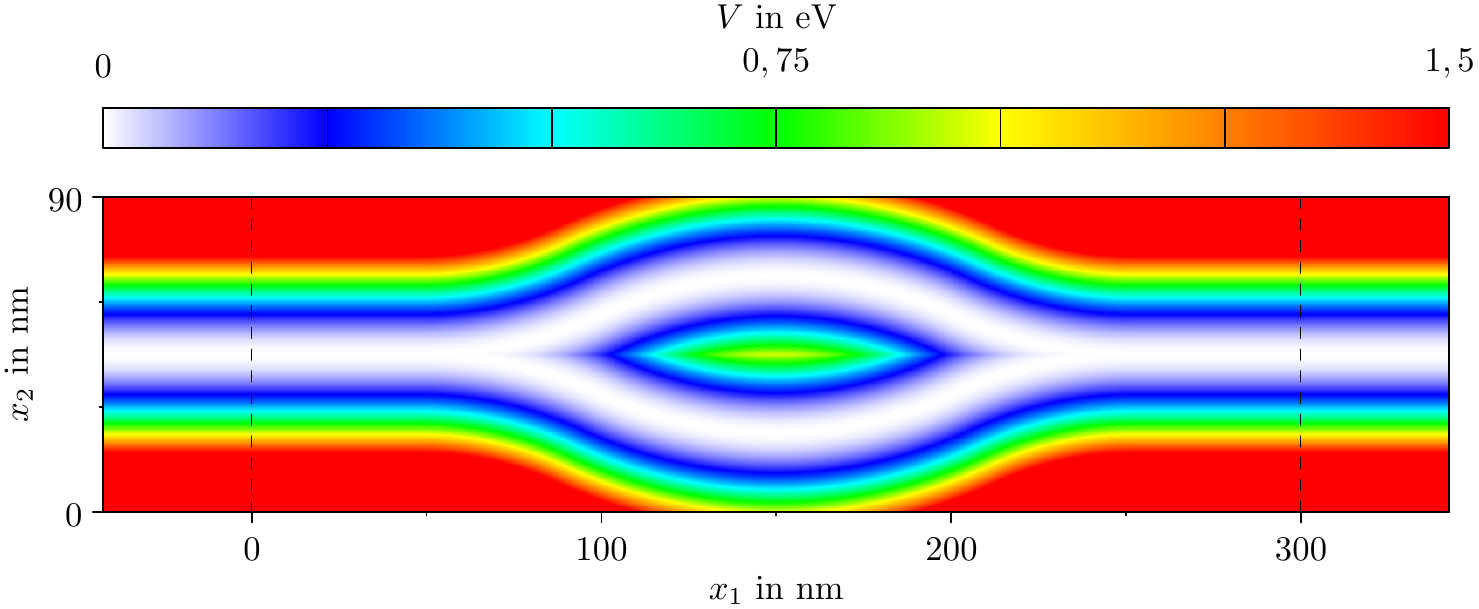}
  	\caption{
  	Potential energy of a ring-shaped quantum waveguide device.
  	The device domain equals $[0,300]\times[0,90]$\,nm$^2$.
	}
  	\label{fig:potential_energy_aharonov_bohm_ring}
  	\end{centering}
\end{figure}

We consider the stationary Schr\"odinger equation
\begin{equation}
  \label{eq:stationary_schroedinger_equation_2d}
  \hat{H} \phi = E \phi, \quad \hat{H} = -\frac{\hbar^2}{2m^\star} \Delta + V,
\end{equation}
on the infinite strip $\Omega = \mathbb{R} \times (0,L_2)$
with homogeneous Dirichlet boundary conditions at 
$x_2=0$ and $x_2=L_2$. For the quantum waveguide simulations,
we assume that the potential energy in the exterior domain depends on the 
transversal coordinate only,
\begin{equation*}
  V(x_1, x_2) =
  \begin{cases}
  \begin{aligned}
  V_\ell(x_2) &\quad \text{for } x_1 \leq 0, \\
  V_r(x_2)    &\quad \text{for } x_1 \geq L_1.
  \end{aligned}
  \end{cases}
\end{equation*}
In general, $V_\ell$ and $V_r$ may be different, but to simplify the notation, 
we assume that $V_\ell(x_2)=V_r(x_2)$ for $x_2 \in [0,L_2]$.
As an example, we refer to the ring-shaped device described by the potential 
energy shown in Figure \ref{fig:potential_energy_aharonov_bohm_ring}.

We further assume that the wave function in the leads, 
\begin{equation}
  \label{eq:decomposition_in_leads_stationary}
  \phi(x_1, x_2) 
  = \sum_{m=0}^{\infty} c^{(m)}(x_1) \, \chi^{(m)}(x_2), \quad x_1 \leq 0, 
	\quad x_1 \geq L_1,
\end{equation}
can be decomposed into transversal waveguide eigenstates,
\begin{equation}
  \label{eq:transversal_waveguide_eigenstates}
  \begin{aligned}
  -\frac{\hbar^2}{2 m^\star} \frac{\partial^2}{\partial x_2^2} \chi^{(m)}(x_2) 
  + V(x_2) \chi^{(m)}(x_2) 
  &= E^{(m)} \chi^{(m)}(x_2), \; \chi^{(m)}(0) = 0, \; \chi^{(m)}(L_2) = 0,\\
  \langle \chi^{(m)}, \chi^{(n)} \rangle 
  &= \int_0^{L_2} \chi^{(m)}(x_2) \, \overline{\chi^{(n)}(x_2)} \, dx_2 = \delta_{m,n},
  \end{aligned}
\end{equation}
where the mode coefficients $c^{(m)}(x_1) = \langle\phi(x_1,\cdot),\chi^{(m)}\rangle$
satisfy the one-dimensional stationary Schr\"odinger equation
\begin{equation}
  \label{eq:stationary_schroedinger_equation_modes}
  -\frac{\hbar^2}{2 m^\star} \frac{\partial^2}{\partial x_1^2} c^{(m)}(x_1)
  + E^{(m)} c^{(m)}(x_1)
  = E c^{(m)}(x_1), \quad m \in \mathbb{N}_0.
\end{equation}

In the following discussion, we consider electrons injected at the left 
terminal traveling to the right.
Let their wave vector be given by $(k,0)^\top$ with $k > 0$.
At the time of the injection, the electrons are assumed to be in the 
ground state with respect to the cross section of the waveguide potential 
in the left lead. Hence, the incoming electrons are represented by
\begin{equation}
  \label{eq:phi_inc_waveguide_stationary}
  \phi^\mathrm{inc}(x_1,x_2) := \exp(i k x_1) \chi^{(0)}(x_2),
\end{equation}
and their total energy amounts to 
$E = \Ekininc + E^{(0)}$, where $\Ekininc = \hbar^2 k^2 / (2 m^\star)$.

\medskip\noindent
{\bf Discrete transparent boundary conditions.}
The symmetric second-order finite-dif\-fer\-ence approximation of the spatial 
derivatives in \eqref{eq:stationary_schroedinger_equation_2d} gives
\begin{equation}
  \begin{aligned}
  \label{eq:discrete_stationary_schroedinger_equation_2d}
  -\frac{\hbar^2}{2 m^\star}\bigg( 
  \frac{\phi_{j_1-1,j_2} - 2 \phi_{j_1,j_2} + \phi_{j_1+1,j_2}}{(\triangle x)^2}
  &+\frac{\phi_{j_1,j_2-1} - 2 \phi_{j_1,j_2} + \phi_{j_1,j_2+1}}{(\triangle x)^2}
  \bigg)\\
  &+ V_{j_1,j_2} \phi_{j_1,j_2} = E \phi_{j_1,j_2}
  \end{aligned}
\end{equation}
on the semi-infinite grid $\Omega_{\triangle x}:= 
\{(j_1 \triangle x, j_2 \triangle x):$ $j_1 \in \mathbb{Z}$, $j_2 = 0,\ldots,J_2\}$,
where $J_1 \triangle x = L_1$, $J_2 \triangle x = L_2$, and 
$\phi_{j_1,0} = \phi_{j_1,J_2} = 0$ for all $j_1 \in \mathbb{Z}$.
In particular, we seek for a solution of 
\eqref{eq:discrete_stationary_schroedinger_equation_2d}
restricted to the grid points of the device domain
\begin{equation*}
  \Omega_\mathrm{DTBC} := \left\{
  (j_1 \triangle x, j_2 \triangle x): j_1 = 0,\ldots,J_1, \; j_2 = 0,\ldots,J_2\right\},
\end{equation*}
where an incoming plane wave is prescribed at the left boundary.
Disregarding open boundary conditions and the incoming plane wave, 
we can state the problem as
\begin{equation}
  \label{eq:discrete_stationary_schroedinger_equation_2d_tensor_product}
  S \phi = 0, \quad \phi_j = \phi_{j_1,j_2}, \quad j = j_1 J_2 + j_2,
  \quad j_1=0,\ldots,J_1, \quad j_2=0,\ldots,J_2.
\end{equation}
Here, $S$ denotes the sparse matrix
$S:= -\hbar^2 / (2 m^\star)\Delta_{x_1,x_2}^\mathrm{2nd}+ \diag(d)$, where
$\Delta_{x_1,x_2}^\mathrm{2nd}:=D_{x_1}^{2,\mathrm{2nd}}\otimes I_{J_2}
+ I_{J_1}\otimes D_{x_2}^{2,\mathrm{2nd}}$.
The components of the vector $d$ are given by 
$d_j=V_j-E,\; j=0,\ldots,(J_1+1)(J_2+1)-1$,
and $I_{J_1}, I_{J_2}$ are unit matrices of dimension $J_1$ and $J_2$, respectively.
Finally, $D_{x_1}^\mathrm{2,2nd}$ and $D_{x_2}^\mathrm{2,2nd}$ are 
finite-difference matrices defined according to \eqref{eq:D_2_2nd} 
with respect to the two spatial directions $x_1$ and $x_2$.

\begin{figure}
	\begin{centering}
  	\includegraphics[width=150mm]{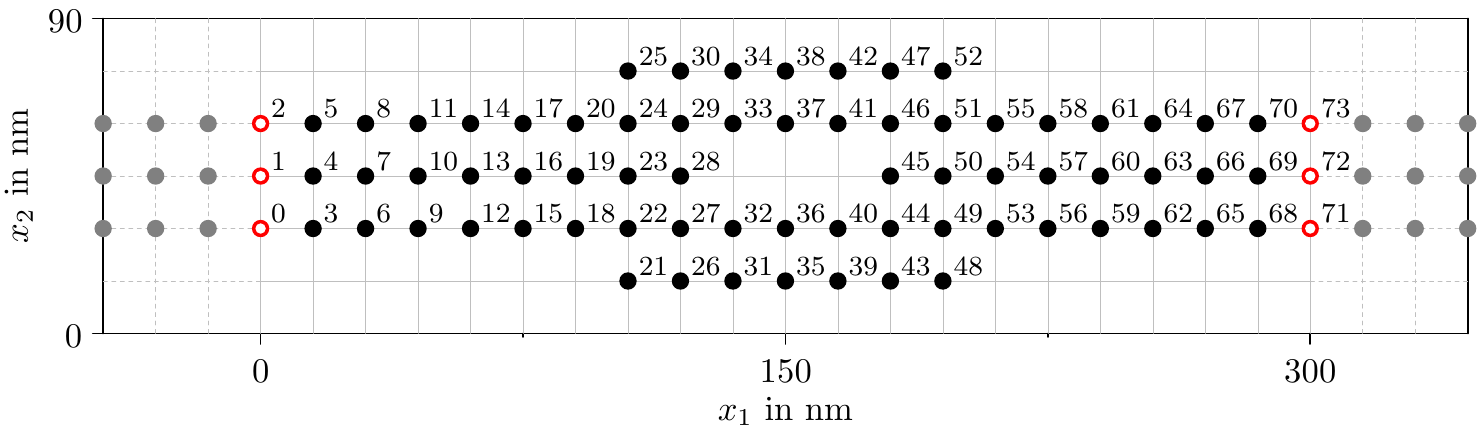}
  	\caption{
  	Reduced mesh of the ring shaped quantum device shown in
  	Figure \ref{fig:potential_energy_aharonov_bohm_ring}.
  	For visualization purposes, $\triangle x$ is chosen extremely large.
  	Grid points at which the wave function is practically zero have been eliminated.
  	For DTBC, the reduced mesh is given by the grid points between 
  	$x_1=0$\,nm and $x_1=300$\,nm.
  	In case of PML, the reduced mesh contains also the grey points.
  	In both cases, we prescribe an incoming plane wave at $x_1=0$\,nm.}
\label{fig:grid_dtbc_pml_2d}
\end{centering}
\end{figure}

In order to realize discrete open boundary conditions at the device terminals, 
we need to replace those finite-difference equations in
\eqref{eq:discrete_stationary_schroedinger_equation_2d_tensor_product}
which correspond to grid points at the left and right boundary of 
$\Omega_\mathrm{DTBC}$.
Furthermore, some of the finite-difference equations 
can be eliminated because $\phi_j$ is zero due to the homogeneous Dirichlet
boundary conditions imposed at the top and the bottom boundary of 
$\Omega_\mathrm{DTBC}$.
In fact, we eliminate even more equations.
Since the wave function $\phi_j$ decreases exponentially within areas where the 
potential energy is greater than the total energy of the electron,
$\phi_j$ is zero (to numerical precision) at some distance from the 
center of the waveguide profile. This allows us to eliminate
the corresponding finite-difference equations.
More specifically, we eliminate all rows $S[j,:]$ and colums 
$S[:,j]$ of $S$ with $j \in \left\{j:V_{j}>750\,\mathrm{meV}\right\}$.
Through this elimination process, we implicitly obtain a reduced mesh
and a new numbering of the remaining free indices.
An illustration is given in Figure \ref{fig:grid_dtbc_pml_2d}.

\begin{figure}[ht]
  	\centering
  	\includegraphics[height=55mm]{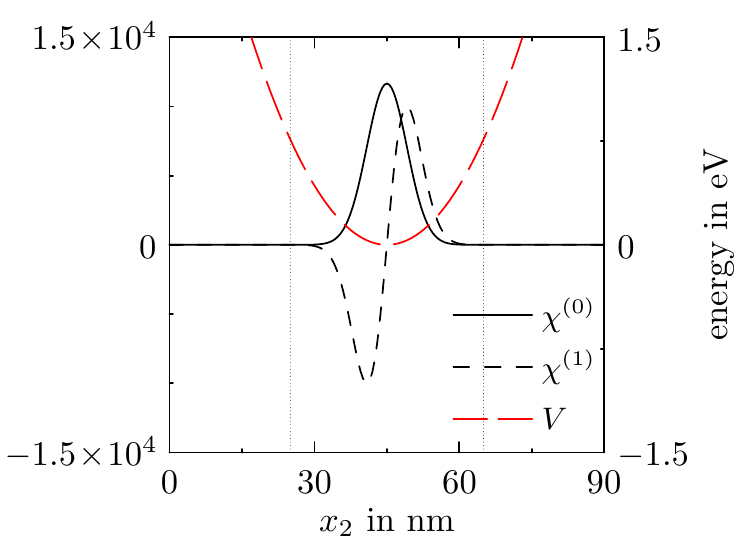}
  	\hspace{15mm}
  	\includegraphics[height=55mm]{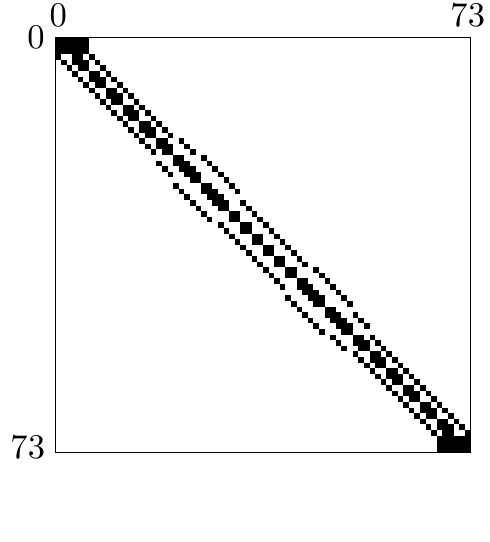}
  	\caption{
  	{\em Left:}
  	Ground state and first excited eigenstate corresponding to the cross-sectional   
	potential energy in the left lead of the quantum waveguide depicted in Figure 
	\ref{fig:potential_energy_aharonov_bohm_ring}. 
  	The vertical dotted lines indicate the boundaries of the reduced mesh.
  	{\em Right:} Sparsity pattern of $S$ 
  	corresponding to the example considered in Figure \ref{fig:grid_dtbc_pml_2d}.}
\label{fig:eigenstates_lhs_and_sparsity_pattern_dtbc}
\end{figure}

We still need to replace the remaining rows of $S$ corresponding to grid points 
at the device contacts (open red points in Figure \ref{fig:grid_dtbc_pml_2d}).
Let us consider the left terminal first.
Since the potential energy in the exterior domain depends solely on the transversal 
coordinate, we temporarily define $V_{j_2}:=V_{j_1,j_2}$ for $j_1\leq 0$ 
and $j_2=0,\ldots,J_2$.
Analogously to the continuous case, the wave function in the lead can be decomposed 
into transversal waveguide eigenstates:
$$
  \phi_{j_1,j_2} = \sum_{m=0}^{M-1} c_{j_1}^{(m)} \chi^{(m)}_{j_2},
  \quad j_1 \leq 0, \quad j_2 = j_{21},\ldots,j_{22}.
$$
The indices $j_{21}$ and $j_{22}$ depend on the elimination process described above.
In the example of Figure \ref{fig:grid_dtbc_pml_2d}, we have $j_{21}=2$ and $j_{22}=4$.
Hence, the number of free indices along the $x_2$-direction in the left lead 
is given by $M=j_{22}-j_{21}+1$.
Further, $\chi^{(m)}$ denotes the $m$-th eigenstate of the discrete eigenvalue problem
\begin{equation*}
  -\frac{\hbar^2}{2 m^\star} 
  \frac{ \chi_{j_2-1}^{(m)} -2 \chi_{j_2}^{(m)} + \chi_{j_2+1}^{(m)}}{(\triangle x)^2} 
  + V_{j_2} \chi_{j_2}^{(m)} = E^{(m)} \chi_{j_2}^{(m)}, 
  \quad j_2=j_{21},\ldots,j_{22}, \ m=0,\ldots,M-1,
\end{equation*}
where we impose homogeneous Dirichlet boundary conditions at 
$x_2 = (j_{21}-1) \triangle x$ and $x_2= (j_{22}+1) \triangle x$.
We further ensure that all eigenstates are orthonormal with respect to the 
scalar product
\begin{equation}
  \label{eq:discrete_l2_scalar_product}
  \langle \chi^{(m)}, \chi^{(n)} \rangle  = 
  \triangle x \sum_{j_2=j_{21}}^{j_{22}} \chi_{j_2}^{(m)} \overline{\chi_{j_2}^{(n)}}.
\end{equation}
The ground state ($m=0$) and the first excited state ($m=1$) are shown
in Figure \ref{fig:eigenstates_lhs_and_sparsity_pattern_dtbc} (left).

In the continuous case, the coefficients $c^{(m)}$ solve 
\eqref{eq:stationary_schroedinger_equation_modes} for all $m \in \mathbb{N}_0$.
The discrete analogue of \eqref{eq:stationary_schroedinger_equation_modes} reads as
\begin{equation}
  \label{eq:stationary_schroedinger_equation_discrete_modes}
  -\frac{\hbar^2}{2 m^\star}
  \frac{c_{j_1-1}^{(m)} - 2 c_{j_1}^{(m)} + c_{j_1+1}^{(m)}}{(\triangle x)^2}
  + E^{(m)} c_{j_1}^{(m)} = E c_{j_1}^{(m)}, \quad m = 0,\ldots,M-1,
\end{equation}
which can be identified with the one-dimensional discrete stationary Schr\"odinger 
equation \eqref{eq:discrete_stationary_schroedinger_equation_1d} if the potential 
energy is substituted by $E^{(m)}$.
Thus, \eqref{eq:stationary_schroedinger_equation_discrete_modes}
admits two solutions of the form 
$c_{j_1}^{(m)} = (\alpha^{(m)})^{j_1}$, where
\begin{align*}
  \alpha^{(m)} 
  &= 1 - \frac{m^\star (E-E^{(m)}) (\triangle x)^2}{\hbar^2} \\
  &\phantom{xx}{}\pm i\sqrt{\frac{2 m^\star (E-E^{(m)})(\triangle x)^2}{\hbar^2}
	-\frac{(m^\star)^2 (E-E^{(m)})^2 (\triangle x)^4}{\hbar^4}}.
\end{align*}
Recall that $m$ denotes the mode index and $m^\star$ the effective mass 
of the electron and that we consider electrons which are injected 
at the left terminal traveling to the right.
At the time of the injection, the electrons are assumed to be in the 
ground state with respect to the waveguide profile.
Thus, we have
$$
  c_{j_1}^{(0)} = A \big( \alpha^{(0)} \big)^{j_1} 
  + B \big( \alpha^{(0)} \big)^{-j_1},  \\
$$
where $A$ and $B$ are the amplitudes of the incoming and the reflected wave 
projected to the ground state $\chi^{(0)}$, respectively.
Writing the above equation for $j_1=0,1$ and eliminating $B$ yields
\begin{equation}
  \label{eq:mode_0_stationary_lhs}
  c_0^{(0)} - \alpha^{(0)} c_1^{(0)} 
  = A \big(1 - \big(\alpha^{(0)}\big)^2\big),
\end{equation}
which corresponds to equation \eqref{eq:dtbc_stationary_1d} with $A=1$.
Here, we use $A=1/\triangle x$ which
gives a reasonable scaling of the final wave function.
In other words, the final wave function 
will be of the same order as the transversal waveguide eigenstates.
In case of the excited modes $(m>0)$, the amplitude of the incoming wave is zero.
Thus, we have
\begin{equation}
  \label{eq:modes_excited_stationary_lhs}
  c_0^{(m)} - \alpha^{(m)} c_1^{(m)} = 0, \quad m = 1,\ldots,M-1.
\end{equation}
Again, the mode coefficients are the projections of the wave function 
onto the transversal waveguide eigenstates. More precisely,
using the scalar product \eqref{eq:discrete_l2_scalar_product} 
and the fact that all eigenstates are real-valued,
we can write the discrete analogue of 
$c^{(m)}$ solving \eqref{eq:stationary_schroedinger_equation_modes} as
\begin{equation}
  \label{eq:discrete_mode_coefficients_as_projection}
  c_{j_1}^{(m)}
  = \triangle x \sum_{j_2=j_{21}}^{j_{22}} \phi_{j_1,j_2} \chi_{j_2}^{(m)}, 
	\quad m=0,\ldots,M-1.
\end{equation}

The remaining free indices corresponding to the grid points at the left boundary
of the device are $0,\ldots,M-1$. Furthermore,
the free indices of the adjacent grid points are $M,\ldots,2M-1$ 
(see Figure \ref{fig:grid_dtbc_pml_2d}). 
Thus, using \eqref{eq:mode_0_stationary_lhs}, \eqref{eq:modes_excited_stationary_lhs}, 
and \eqref{eq:discrete_mode_coefficients_as_projection}, DTBC at $x_1=0$ become
\begin{equation}
  \label{dtbc_2d_left_terminal}
  \bigg( \triangle x \sum_{j=0}^{M-1} \phi_j \chi_j^{(m)} \bigg)
  -\alpha^{(m)} \bigg( \triangle x \sum_{j=M}^{2M-1} \phi_j \chi_j^{(m)} \bigg)
  = \frac{1}{\triangle x}\big(1 - \big(\alpha^{(0)}\big)^2\big) \delta_{m,0},
\end{equation}
for all $m=0,\ldots,M-1$.
Each of these $M$ equations is used to replace one equation corresponding 
to a grid point at the left boundary of the reduced mesh.
In other words, we replace the first $M$ rows of $S$ according to the left-hand side of
\eqref{dtbc_2d_left_terminal}.

DTBC at $x_1=L_1$ follow analogously.
Since we do not prescribe an incoming wave at the right contact, the inhomogeneity 
does not show up, leading to
\begin{equation}
  \label{dtbc_2d_right_terminal}
  \bigg( \triangle x \sum_{j=N-M}^{N-1} \phi_j \chi_j^{(m)} \bigg)
  -\alpha^{(m)} \bigg( \triangle x \sum_{j=N-2M}^{N-M-1} \phi_j \chi_j^{(m)} \bigg) 
  = 0, \quad m=0,\ldots,M-1.
\end{equation}
The last $M$ rows of $S$ are replaced according to \eqref{dtbc_2d_right_terminal}
which finally yields
\begin{equation}
  \label{eq:quantum_waveguide_stationary_dtbc_final_equations}
  S \phi = b, \quad \phi=(\phi_0, \ldots, \phi_{N-1})^\top, 
  \quad b = (b_0, \ldots, b_{N-1})^\top,
\end{equation}
where $b_j = (1/\triangle x)(1 - (\alpha^{(0)})^2) \delta_{j,0}$.

In some numerical experiments we noticed that the condition number of $S$ is 
quite large. This problem can be easily overcome by scaling all equations
(apart from the ones given in \eqref{dtbc_2d_left_terminal}
and \eqref{dtbc_2d_right_terminal}) with $1/E$.
The sparsity pattern of $S$ is depicted in Figure
\ref{fig:eigenstates_lhs_and_sparsity_pattern_dtbc}.
Obviously, $S$ contains two dense submatrices which are a direct consequence of the
DTBC \eqref{dtbc_2d_left_terminal} and \eqref{dtbc_2d_right_terminal}.
For that reason, the symmetry of $S$ is lost and consequently, many iterative 
methods cannot be applied to solve
\eqref{eq:quantum_waveguide_stationary_dtbc_final_equations}.
However, the simulations considered below can still be handled by direct solvers
and we are not affected by the loss of symmetry of $S$.

\medskip\noindent
{\bf Perfectly Matched Layers.}
Similarly to the one-dimensional case, we replace the Laplacian in the 
stationary Schr\"odinger equation by the Laplace-PML operator
\begin{equation}
  \label{eq:laplace_pml_operator}
  \frac{\partial^2}{\partial x_1^2} + \frac{\partial^2}{\partial x_2^2}
  \rightarrow  c(x_1)\frac{\partial}{\partial x_1}c(x_1)
  \frac{\partial}{\partial x_1}+\frac{\partial^2}{\partial x_2^2},
\end{equation}
which yields the stationary Schr\"odinger-PML equation
\begin{equation}
  \label{eq:stationary_schroedinger_pml_equation_2d}
  -\frac{\hbar^2}{2 m^\star}
  \left(c(x_1)\frac{\partial}{\partial x_1}c(x_1)
  \frac{\partial}{\partial x_1}
  +\frac{\partial^2}{\partial x_2^2} \right)
  \phi(x_1,x_2) + V(x_1, x_2) \phi(x_1,x_2) 
  = E \phi(x_1,x_2).
\end{equation}
We choose the same function $c(x)$ as in Section \ref{subsec:scattering_states_1d}.

According to \eqref{eq:D_tilde_1d_2nd}, a second-order finite-difference 
discretization of  \eqref{eq:stationary_schroedinger_pml_equation_2d} is given by
\begin{equation*}
  \begin{aligned}
  -\frac{\hbar^2}{2 m^\star}
  &\bigg(
  c(x_{j_1}) c'(x_{j_1}) \frac{-\phi_{j_1-1,j_2} + \phi_{j_1+1,j_2}}{2 \triangle x}
  + c^2(x_{j_1}) \frac{\phi_{j_1-1,j_2}-2 \phi_{j_1,j_2} 
	+ \phi_{j_1+1,j_2}}{(\triangle x)^2} \\
  &+\frac{\phi_{j_1,j_2-1} - 2 \phi_{j_1,j_2} + \phi_{j_1,j_2+1}}{(\triangle x)^2}
  \bigg) + V_{j_1,j_2} \phi_{j_1,j_2}
  = E \phi_{j_1,j_2}, \;\; j_1 \in \mathbb{Z}, \;\; j_2 = 0,\ldots,J_2.
  \end{aligned}
\end{equation*}
In fact, we seek for a solution restricted to the computational domain
\begin{equation*}
  \Omega_\mathrm{PML} = \left\{(X_{j_1}^\mathrm{PML}, j_2 \triangle x):
  j_1 = 0,\ldots,J_1^\mathrm{PML}, \; j_2 = 0,\ldots,J_2\right\},
\end{equation*}
where $X_{j_1}^\mathrm{PML}$ denotes the $j_1$-th grid point of the one-dimensional 
grid $X_\mathrm{PML}$ defined in \eqref{eq:X_pml}.
Without taking into account a possible incoming plane wave, we can state the problem 
in the following form:
\begin{equation*}
  S_\mathrm{PML} \phi = 0, \quad \phi_j = \phi_{j_1,j_2}, \quad j = j_1 J_2 + j_2, 
  \quad j_1=0,\ldots,J_1^\mathrm{PML}. \quad j_2=0,\ldots,J_2. 
\end{equation*}
where $S_\mathrm{PML}:=-\hbar^2 / (2 m^\star)\tilde{\Delta}_{x_1,x_2} + \diag(d)$,
and $\tilde{\Delta}_{x_1,x_2}$ is one of the sparse matrices
\begin{subequations}
  \label{eq:laplacians_pml}
  \begin{align}
  \begin{split}
  \label{eq:laplacian_pml_2nd_order}
  \tilde{\Delta}_{x_1,x_2}^\mathrm{2nd}
	:=\left(\tilde{D}_{x_1}^{2,\mathrm{2nd}}\otimes I_{J_2}
  + I_{J_1^\mathrm{PML}}\otimes D_{x_2}^{2, \mathrm{2nd}}\right),
  \end{split} \\
  \begin{split}
  \label{eq:laplacian_pml_4th_order}
  \tilde{\Delta}_{x_1,x_2}^\mathrm{4th}
	:=\left(\tilde{D}_{x_1}^{2,\mathrm{4th}}\otimes I_{J_2}
  + I_{J_1^\mathrm{PML}}\otimes D_{x_2}^{2, \mathrm{4th}} \right),
  \end{split} \\
  \begin{split}
  \label{eq:laplacian_pml_6th_order}
  \tilde{\Delta}_{x_1,x_2}^\mathrm{6th}
	:=\left(\tilde{D}_{x_1}^{2,\mathrm{6th}}\otimes I_{J_2}
	+ I_{J_1^\mathrm{PML}} \otimes D_{x_2}^{2, \mathrm{6th}}\right),
  \end{split}
  \end{align}
\end{subequations}
corresponding to second-, fourth-, and sixth-order discretizations of the 
Laplace-PML operator \eqref{eq:laplace_pml_operator}.
The finite-difference matrices corresponding to the $x_1$-direction are 
defined via \eqref{eq:D_tilde_1d} and the finite-difference matrices 
acting on the $x_2$-direction are given in \eqref{eq:D_2_2nd}-\eqref{eq:D_2_6th}.
Finally, the vector $d$ is defined via $d_j=V_j-E$, 
$j=0,\ldots,(J_1^\mathrm{PML}+1)(J_2+1)-1$.

Before we realize an incoming plane wave at the left terminal, we eliminate
all finite-difference equations corresponding to grid points where the wave 
function is supposed to be zero. Like in the case of DTBC, we eliminate
all rows $S_\mathrm{PML}[j,:]$ and colums $S_\mathrm{PML}[:,j]$ of 
$S_\mathrm{PML}$ with $j \in \left\{j:V_{j}>750\,\mathrm{meV}\right\}$.
Let $N_\mathrm{PML} \times N_\mathrm{PML}$ denote the new size of $S_\mathrm{PML}$.
The elimination process implicitly yields a reduced mesh which extends the reduced 
mesh considered in the case of DTBC (see Figure \ref{fig:grid_dtbc_pml_2d}). 
The wave function in the left contact is a superposition of an incoming and 
a reflected wave,
$\phi_{j_1,j_2} = \phi_{j_1,j_2}^\mathrm{inc} + \phi_{j_1,j_2}^\mathrm{refl}$.
The incoming wave 
\begin{equation}
  \label{eq:inc_wave_pml_stationary_2d}
  \phi_{j_1,j_2}^\mathrm{inc} = \chi_{j_2}^{(0)} \exp(i k j_1 \triangle x),
  \quad j_1 \leq 0, \quad j_2=j_{21},\ldots,j_{22},
\end{equation}
is a discrete representation of \eqref{eq:phi_inc_waveguide_stationary}, where 
$\chi^{(0)}$ is a solution of
$$
  -\frac{\hbar^2}{2 m^\star} D_{x_2}^2 \chi_{j_2}^{(0)} + V_{j_2} 
	\chi_{j_2}^{(0)} = E^{(0)}, 
  \quad D_{x_2}^2 \in \big\{D_{x_2}^{2,\mathrm{2nd}},D_{x_2}^{2,\mathrm{4th}}, 
	D_{x_2}^{2,\mathrm{6th}}\big\}, \quad j_2=j_{21},\ldots,j_{22}.
$$
Like in the case of DTBC, we impose homogeneous Dirichlet boundary conditions at 
$x_2=(j_{21}-1) \triangle x$ and $x_2=(j_{22}+1) \triangle x$.
The wave number $k$ is related to the kinetic energy $\Ekininc$ according to the 
discrete $E$--$k$--relation \eqref{eq:discrete_E_k_relation}.
In case of the higher-order methods, we simply use the continuous relation
\eqref{eq:continuous_E_k_relation}.
Let again $M$ denote the number of free indices along the $x_2$-direction 
in the left lead.
Moreover, let $j_0,\ldots,j_0+M-1$ denote the free indices of the grid points at 
$x_1=0$ with respect to the reduced mesh 
($j_0=9$ in Figure \ref{fig:grid_dtbc_pml_2d}).
With these notations, an incoming plane wave is realized in the same way as in the 
one-dimensional case. For the second-order discretization, we find that
\begin{equation*}
  S_\mathrm{PML} \phi = b, \quad \phi=(\phi_0, \ldots, \phi_{N_\mathrm{PML}-1})^\top, 
  \quad b = (b_0, \ldots, b_{N_\mathrm{PML}-1})^\top,
\end{equation*}
where
$$
  b_j =
  \begin{cases}
  \begin{aligned}
  -\frac{\hbar^2}{2 m^\star (\triangle x)^2} 
	\phi_{j+M}^\mathrm{inc} &\quad  \text{for } j = j_0-M,\ldots,j_0-1, \\
  +\frac{\hbar^2}{2 m^\star (\triangle x)^2} 
	\phi_{j-M}^\mathrm{inc} &\quad  \text{for } j = j_0,\ldots,j_0+M-1, \\
  0 &\quad \textrm{else}.
  \end{aligned}
  \end{cases}
$$
In case of the higher-order discretizations, $b$ needs to be adapted accordingly.

Since $S_\mathrm{PML}$ does not contain dense submatrices,
it can be assembled more easily than the corresponding matrix in 
\eqref{eq:quantum_waveguide_stationary_dtbc_final_equations}.
However, depending on the ratio between the number of grid points needed to 
describe the device domain and the number of grid points needed to realize the PML,
the size of $S_\mathrm{PML}$ may be significantly larger than that of 
$S$ in \eqref{eq:quantum_waveguide_stationary_dtbc_final_equations}.
Moreover, $S_\mathrm{PML}$ is not Hermitian. 
Thus, like in the case of DTBC, many iterative methods cannot be applied.

\medskip\noindent
{\bf Simulations.} 
Before we turn our attention to the ring-shaped quantum device introduced above,
let us first consider a straight waveguide with parabolic cross section, i.e.,
\begin{equation}
  \label{eq:waveguide_parabolic}
  V(x_1,x_2) = \frac12 m^\star \omega_\ast^2 \left(x_2-\frac{L_2}{2}\right)^2,
  \quad \omega_\ast = 0.5 \times 10^{14}\mathrm{\,s}^{-1},
  \ x_1 \in \mathbb{R}, \ x_2 \in [0,L_2].
\end{equation}
The cross-sectional eigenstates and eigenvalues read as
\begin{equation}
  \label{eq:parabolic_cross_sectional_eigenstates}
  \begin{aligned}
  \chi^{(n)}(x_2) =& \left( \frac{m^\star \omega_\ast}{\pi \hbar} \right)^{1/4}
  \frac{1}{\sqrt{2^n n !}} 
	H_n\left( \sqrt{\frac{m^\star \omega_\ast}{\hbar}} 
	\left(x_2-\frac{L_2}{2}\right)\right)\\
  &\times\exp\left(-\frac{1}{2} \frac{m^\star \omega_\ast}{\hbar} 
	\left(x_2-\frac{L_2}{2}\right)^2\right), \quad 
	E^{(n)} = \hbar \omega_\ast \left( n + \frac{1}{2}\right), \quad n \in \mathbb{N}_0,
  \end{aligned}
\end{equation}
where $H_n$ denote the Hermite polynomials.

\begin{figure}
	\centering
  	\includegraphics[width=75mm]{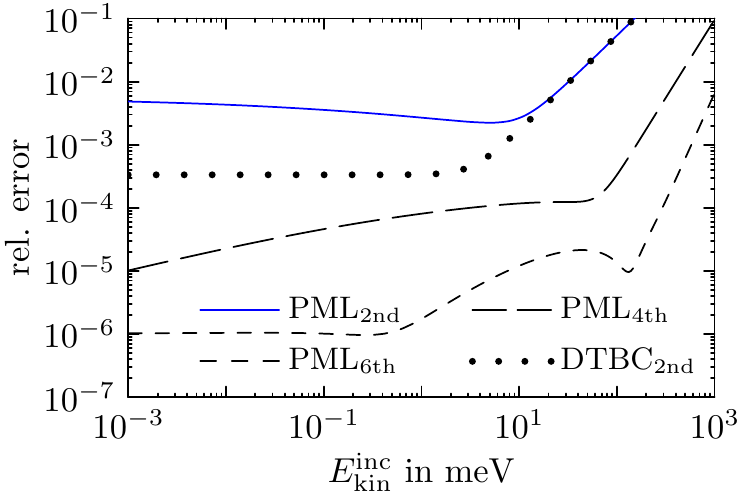}
  	\hspace{2.5mm}
  	\includegraphics[width=75mm]{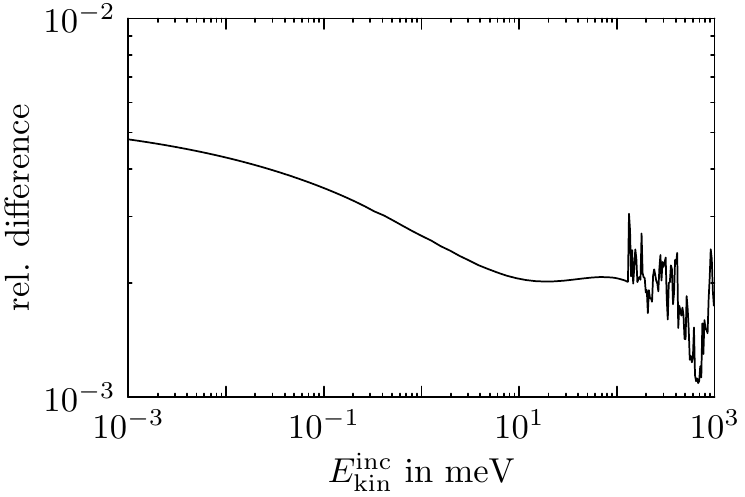}
  	\caption{
  	{\em Left:}
  	Relative errors as a function of $\Ekininc$ in a numerical scattering 
	state experiment. The electrons are injected at the left contact of a straight 
	waveguide with parabolic cross section.
  	{\em Right:}
  	Relative differences between the solutions obtained by $\mathrm{PML}_\mathrm{2nd}$ 
	and $\mathrm{DTBC}_\mathrm{2nd}$ as a function of $\Ekininc$.
  	The electrons are injected at the left terminal of the ring-shaped waveguide 
	depicted in Figure \ref{fig:potential_energy_aharonov_bohm_ring}.}
  \label{fig:rel_errors_differences_scattering_states_2d}
\end{figure}

In this trivial example, the solution of the scattering state problem  
is given by \eqref{eq:phi_inc_waveguide_stationary}, where $\chi^{(0)}$ is 
substituted according to \eqref{eq:parabolic_cross_sectional_eigenstates}.
We compare the exact solution with the results of the different numerical solvers
using the mesh size $\triangle x = 0.5$\,nm.
To this end, we set $[0,L_1]\times[0,L_2] = [0,120]\times[0,60]\,\mathrm{nm}^2$.
The relative errors are depicted in Figure 
\ref{fig:rel_errors_differences_scattering_states_2d} (left)
for the energy range $\Ekininc \in [10^{-3},10^3]$\,meV.
As expected, the results resemble the results from the corresponding one-dimensional 
simulation (Figure \ref{fig:rel_errors_scattering_states_1d} right).
However, while in the one-dimensional simulation, the errors of 
$\mathrm{DTBC}_\mathrm{2nd}$ are decreasing continuously for decreasing $\Ekininc$, 
they are now bounded from below.
This lower bound results from the finite numerical resolution of the transversal 
waveguide eigenstate $\chi^{(0)}$.

\begin{figure}[ht]
	\begin{centering}
  	\includegraphics[width=150mm]{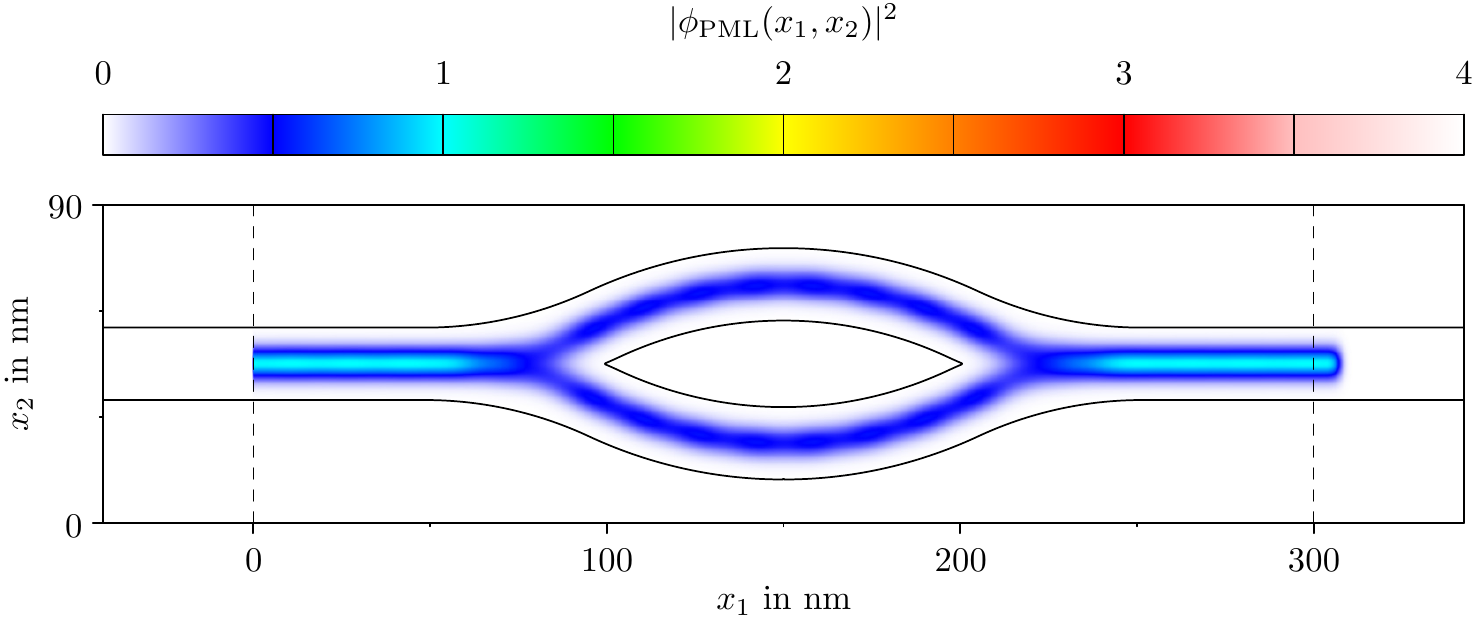}
  	\caption{
  	Scattering state in the ring-shaped quantum waveguide shown in Figure
  	\ref{fig:potential_energy_aharonov_bohm_ring}.
  	The kinetic energy of the electrons injected at the left contact is 
	$\Ekininc=21.5$\,meV. The black solid line indicates an isoline of the 
	potential energy at $200$\,meV.}
  	\label{fig:scattering_state_ab_ring_E_kin_inc_21_5_meV}
  	\end{centering}
\end{figure}

We now consider the ring-shaped quantum device depicted in Figure 
\ref{fig:potential_energy_aharonov_bohm_ring}.
A scattering state solution according to $\mathrm{PML}_\mathrm{2nd}$ is shown in 
Figure \ref{fig:scattering_state_ab_ring_E_kin_inc_21_5_meV} for 
$\Ekininc=21.5$\,meV. Here and in all subsequent figures,
the wave function is scaled by the maximum of the transversal waveguide 
eigenstate $\chi^{(0)}$.
As can be seen, the electrons are transmitted almost perfectly through the device. 
Hence, the reflected wave in the left terminal is practically zero.

Next, we compute scattering state solutions for the energy range
$\Ekininc \in [10^{-3},10^3]$\,meV. Since exact solutions are not available and
our solver using DTBC is restricted to second-order accuracy,
we only compare the solutions corresponding to
$\mathrm{DTBC}_\mathrm{2nd}$ and $\mathrm{PML}_\mathrm{2nd}$. 
Their relative differences at a spatial resolution of $\triangle x=0.5$\,nm are 
given in Figure \ref{fig:rel_errors_differences_scattering_states_2d} (right).
For energies up to approximately $100$\,meV, the results
are consistent with the results shown in 
Figure \ref{fig:rel_errors_differences_scattering_states_2d} (left).
For larger energies, the small wavelength of the wave function can hardly 
be resolved with second-order methods.
Nonetheless, $\mathrm{DTBC}_\mathrm{2nd}$ and $\mathrm{PML}_\mathrm{2nd}$ yield 
approximately the same results, even though the relative differences vary 
considerably. 

\begin{figure}
	\centering
  	\includegraphics[width=75mm]{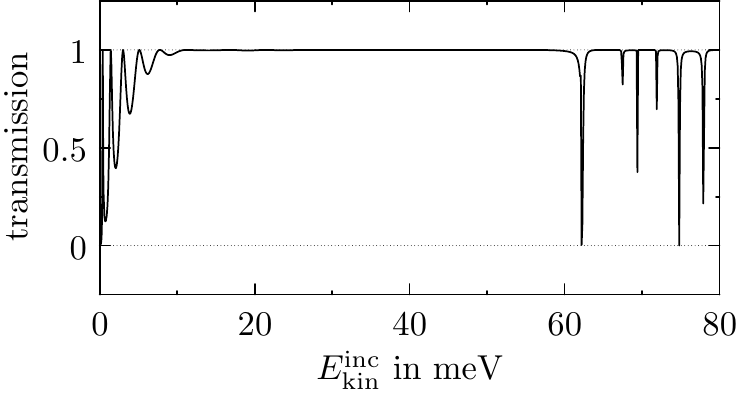}
  	\caption{
		Transmission probability as a function of $\Ekininc$ for electrons 
  		injected at the left terminal of the ring-shaped device shown in 
  		Figure \ref{fig:potential_energy_aharonov_bohm_ring}.
	}
  	\label{fig:transmissions_ab_ring1}
\end{figure}

Finally, we compute the transmission probability as a 
function of the kinetic energy of the incident electrons.
The transmission probability is defined as the ratio between the transmitted and 
the incident probability current density \cite{Dav98}.
In the given situation (provided $\Ekininc$ is not too large) we have
\begin{equation}
  \label{eq:transmission_probability}
  \frac{j_\mathrm{trans}}{j_\mathrm{inc}}
  = \big|\langle \phi(L_1,\cdot), \chi^{(0)} \rangle\big|^2,
\end{equation}
which is depicted in Figure \ref{fig:transmissions_ab_ring1}.
 
\subsection{Transient scattering states in quantum waveguides}

We consider the time-dependent Schr\"odinger equation
\begin{equation}
  \label{eq:tdse_2d}
  i \hbar \partial_t \psi = \hat{H} \psi, \quad 
	\hat{H} = -\frac{\hbar^2}{2 m^\star} \Delta + V
\end{equation}
on the infinite strip $\Omega = \mathbb{R} \times (0,L_2)$.
Like in the stationary case, we prescribe homogeneous Dirichlet boundary
conditions at $x_2=0$ and $x_2=L_2$.
Moreover, we assume that $V(x_1,x_2,t) = V(x_2)$ for 
$x_1\le 0$ or $x_1\ge L_1$ and $x_2\in[0,L_2]$.
This condition is satisfied for the potential in Figure 
\ref{fig:potential_energy_aharonov_bohm_ring}.
In many applications, the potentials in the left and right lead 
are not necessarily the same.
Furthermore, they may depend on time via an applied voltage, e.g.\
$V(x_1,x_2,t) = V(x_2) - e U(t)$.
These extensions can be easily included in the following discussion
using the results from Section \ref{sec:one_dimensional_simulations}.

In the transient scattering state experiment, the initial wave function is given 
by a scattering state solution of the stationary Sch\"odinger equation
for the potential energy $V(x_1,x_2,0)$ and the total energy $E=\Ekininc+E^{(0)}$,
where $\Ekininc=\hbar^2 k^2 / (2 m^\star)$.

Using time-dependent mode coefficients, the wave function in the leads can still 
be decomposed into transversal waveguide eigenstates 
\eqref{eq:transversal_waveguide_eigenstates},
\begin{equation}
  \label{eq:decomp_leads_transient}
  \psi(x_1, x_2, t) = \sum_{m=0}^{\infty} d^{(m)}(x_1,t) \, \chi^{(m)}(x_2), 
	\quad x_1   \leq 0, \quad x_1 \geq L_1.
\end{equation}
Substituting \eqref{eq:decomp_leads_transient} into \eqref{eq:tdse_2d} shows that
each coefficient satisfies the one-dimensional time-dependent Schr\"odinger equation
\begin{equation}
  \label{eq:tdse_modes}
  i \hbar \frac{\partial}{\partial t} d^{(m)}(x_1,t) = 
  -\frac{\hbar^2}{2 m^\star} \frac{\partial^2}{\partial x_1^2} d^{(m)}(x_1,t) 
	+ E^{(m)} d^{(m)}(x_1,t), \quad t \geq 0, \quad m\in\mathbb{N}_0.
\end{equation}

As an eigenstate of the Schr\"odinger equation, the scattering state evolves in time 
according to $\exp(-i \omega t)$ with $\omega = E / \hbar$.
Hence, the time-evolution of the scattering state in the leads is given by
\begin{equation*}
\sum_{m=0}^{\infty} \exp(-i E t / \hbar)\, c^{(m)}(x_1)\, \chi^{(m)}(x_2), 
\quad x_1 \leq 0, \quad x_1 \geq L_1,
\end{equation*}
where $e^{(m)}(x_1,t) := \exp(-i E t / \hbar) c^{(m)}(x_1)$
solves \eqref{eq:tdse_modes} as well.
Hence, $\exp(i E^{(m)} t / \hbar) d^{(m)}$ and 
$\exp(i E^{(m)} t / \hbar) e^{(m)}$
solve the free time-dependent one-dimensional Schr\"o\-dinger equation
and therefore transparent boundary conditions at $x_1=0$ and $x_1=L_1$ can 
be derived by the application of \eqref{eq:tbc_zero_exterior_potential_1d} to
\begin{equation}
  \label{eq:d_minus_e}
  \exp(i E^{(m)} t / \hbar) d^{(m)}-\exp(i E^{(m)} t / \hbar) e^{(m)}
\end{equation}
for each $m \in \mathbb{N}_0$.

\medskip\noindent
{\bf Discrete transparent boundary conditions.}
Using the same grid $\Omega_\mathrm{DTBC}$ as in the stationary case, 
we formulate the Crank-Nicolson scheme for the two-dimensional 
time-dependent Schr\"odinger equation \eqref{eq:tdse_2d} as follows:
\begin{equation}
\label{eq:discrete_tse_2d}
  P \psi^{(n+1)} = Q \psi^{(n)}, \quad \psi_j^{(n)} = \psi_{j_1,j_2}^{(n)},
  \quad j = j_1 J_2 + j_2,
\end{equation}
where $j_1=0,\ldots,J_1$, $j_2=0,\ldots,J_2$, $n \in \mathbb{N}_0$,
$P$ and $Q$ are the sparse matrices
$$
\begin{aligned}
  P:= I-\frac{i \triangle t \hbar}{4 m^\star}\Delta_{x_1,x_2}^\mathrm{2nd}
  + \frac{i \triangle t}{2 \hbar} 
  \diag \Big(
  \big(V_0^{(n+1/2)},\ldots,V_{J_1 J_2-1}^{(n+1/2)}\big)^\top\Big), \\
  Q:= I+\frac{i \triangle t \hbar}{4 m^\star}\Delta_{x_1,x_2}^\mathrm{2nd}
  -\frac{i \triangle t}{2 \hbar}
  \diag  \Big(
  \big(V_0^{(n+1/2)},\ldots,V_{J_1 J_2-1}^{(n+1/2)}\big)^\top\Big).
\end{aligned}
$$
As outlined in Section \ref{subsec:scattering_states_waveguide},
we eliminate all finite-difference equations corresponding to grid points 
at which the wave function is zero.
Needless to say, we eliminate the same rows and columns of $P$ and $Q$ 
which have been eliminated in $S$ (see Section 
\ref{subsec:scattering_states_waveguide}).
Using the same notations as in Section \ref{subsec:scattering_states_waveguide},
the wave function in the leads of the reduced mesh reads as
$$
  \psi_{j_1,j_2}^{(n)} 
	= \sum_{m=0}^{M-1} d_{j_1}^{(m,n)} \chi^{(m)}_{j_2}, 
	\quad j_1 \leq 0, \quad j_1 \geq J_1, \quad j_2 = j_{21},\ldots,j_{22}.
$$

In order to derive DTBC at $x_1=0$ and $x_1=L_1$, we employ the same strategy 
as in the continuous case considered above.
The discrete analogue of \eqref{eq:d_minus_e} becomes
\begin{equation}
  \label{eq:epsilon_minus_gamma}
  \epsilon^{(m,n)} d_{j_1}^{(m,n)} - \gamma^{(m,n)} c_{j_1}^{(m)},
\end{equation}
where the discrete gauge change terms
\begin{equation*}
  \begin{aligned}
  \epsilon^{(m,n)}
  &= \exp \big(2 i n\arctan ( \triangle t E^{(m)} / (2 \hbar) ) \big)
  \approx \exp\big(i E^{(m)} t / \hbar\big),\\
	\gamma^{(m,n)} 
	&= \exp (2 i n ( \arctan (\triangle t E^{(m)} / (2 \hbar) ) 
  - \arctan ( \triangle t E / ( 2 \hbar ) ) ))\\
  &\approx \exp(i E^{(m)} t / \hbar) \exp(-i E t / \hbar),
  \quad m =0,\ldots,M-1, \quad n \in \mathbb{N}_0,
  \end{aligned}
\end{equation*}
are slight modifications of \eqref{eq:gauge_change_gamma_and_epsilon}.
Applying \eqref{eq:hom_dtbc_1d_lhs} and \eqref{eq:hom_dtbc_1d_rhs} to
\eqref{eq:epsilon_minus_gamma} yields the DTBC
\begin{subequations}
  \label{eq:non_hom_dtbc_modes}
  \begin{align}
  \begin{split}
  &\begin{aligned} \label{eq:non_hom_dtbc_modes_a}
  \epsilon^{(m,n+1)} d_1^{(m,n+1)}
  &-s^{(0)}\epsilon^{(m,n+1)} d_0^{(m,n+1)} 
  = \sum_{\ell=1}^{n} s^{(n+1-\ell)} 
  \big(\epsilon^{(m,\ell)} d_0^{(m,\ell)} - \gamma^{(m,\ell)} c_0^{(m)}\big)\\
  &-\big(\epsilon^{(m,n)} d_1^{(m,n)} - \gamma^{(m,n)} c_1^{(m)}\big)
	+\gamma^{(m,n+1)} c_1^{(m)} - s^{(0)} \gamma^{(m,n+1)} c_0^{(m)},
  \end{aligned}
  \end{split} \\
  \begin{split}
  &\begin{aligned} \label{eq:non_hom_dtbc_modes_b}
  \epsilon^{(m,n+1)} d_{J_1-1}^{(m,n+1)}
  &-s^{(0)}\epsilon^{(m,n+1)} d_{J_1}^{(m,n+1)}
  = \sum_{\ell=1}^{n} s^{(n+1-\ell)} 
  \big(\epsilon^{(m,\ell)} d_{J_1}^{(m,\ell)} - \gamma^{(m,\ell)} c_{J_1}^{(m)}\big) \\
  &- \big(\epsilon^{(m,n)} d_{J_1-1}^{(m,n)} - \gamma^{(m,n)} c_{J_1-1}^{(m)}\big)
  + \gamma^{(m,n+1)} c_{J_1-1}^{(m)} - s^{(0)} \gamma^{(m,n+1)} c_{J_1}^{(m)}
  \end{aligned}
  \end{split}
  \end{align}
\end{subequations}
at the left and right contact, respectively.
The mode coefficients appearing on the left-hand side of \eqref{eq:non_hom_dtbc_modes}
are implicitly given by the projection of the new wave function onto the transversal 
waveguide eigenstates
\begin{align*}
  d_0^{(m,n+1)} &= \triangle x \sum_{j=0}^{M-1} \psi_j^{(n+1)} \chi_j^{(m)}, &
  d_1^{(m,n+1)} &= \triangle x \sum_{j=M}^{2M-1} \psi_j^{(n+1)} \chi_j^{(m)}, \\
  d_{J_1-1}^{(m,n+1)} &= \triangle x \sum_{j=N-2M}^{N-M-1} 
	\psi_j^{(n+1)} \chi_j^{(m)}, &
  d_{J_1}^{(m,n+1)} &= \triangle x \sum_{j=N-M}^{N-1} \psi_j^{(n+1)} \chi_j^{(m)}.
\end{align*}

Finally, all remaining equations in \eqref{eq:discrete_tse_2d}, which correspond to 
grid points at the left and right device contacts, are replaced.
In the particular example of the reduced mesh shown in Figure 
\ref{fig:grid_dtbc_pml_2d},
the first $M$ equations are replaced according to \eqref{eq:non_hom_dtbc_modes_a}
and the last $M$ equations are replaced according to \eqref{eq:non_hom_dtbc_modes_b}.
The left-hand sides of \eqref{eq:non_hom_dtbc_modes} cause dense submatrices in $P$ 
which need to be updated in each time step.
All quantities appearing on the right-hand side of \eqref{eq:non_hom_dtbc_modes} are 
already known (at the $n$-th time step) or can be easily computed.
However, to be able to perform the next time step,
we need to store all values of $d_0^{(m,\ell)}$ and $d_{J_1}^{(m,\ell)}$ for
$m=0,\ldots,M-1$ and $\ell=1,\ldots,n$.
Hence, the storage requirements for the DTBC
are in $O(M n_\star)$, where $n_\star$ denotes the total number of time steps.
Furthermore, the computational time required for a single evaluation of the 
discrete convolutions in \eqref{eq:non_hom_dtbc_modes} increases linearly 
with $n$ and thus, the total time complexity is in $O(M n_\star^2)$.

\medskip\noindent
{\bf Perfectly Matched Layers.}
Substituting the Laplacian in the time-dependent Schr\"o\-dinger equation 
with the Laplace-PML operator \eqref{eq:laplace_pml_operator} gives the 
time-dependent Schr\"odinger-PML equation for $\psi=\psi(x_1,x_2,t)$:
\begin{equation*}
  i \hbar \frac{\partial\psi}{\partial t}
  = -\frac{\hbar^2}{2 m^\star}\left(c(x_1)\frac{\partial}{\partial x_1}
  c(x_1)\frac{\partial}{\partial x_1}+\frac{\partial^2}{\partial x_2^2}\right)
  \psi + V(x_1, x_2, t) \psi.
\end{equation*}
The semi-discretized problem on $\Omega_\mathrm{PML}$ reads as
\begin{equation*}
  \frac{d}{d t} \psi_j(t) = 
  i \frac{\hbar}{2 m^\star} \tilde{\Delta}_{x_1,x_2} \psi_j(t) 
  - \frac{i}{\hbar} V(t) \psi_j(t), \quad \psi_j(t) = \psi_{j_1,j_2}(t),
\end{equation*}
where $j = j_1 J_2 + j_2$, $j_1=0,\ldots,J_1^\mathrm{PML}$, $j_2=0,\ldots,J_2$, and
$\tilde{\Delta}_{x_1,x_2}$ denotes one of the sparse matrices defined in 
\eqref{eq:laplacians_pml}.
For the time-integration method, we employ the Crank-Nicolson scheme
\eqref{eq:Crank_Nicolson_time_integration}
or the classical Runge-Kutta method.
Like in the stationary case, we eliminate all equations
corresponding to grid points at which the wave function is assumed to be zero.
The incoming plane wave \eqref{eq:inc_wave_pml_stationary_2d} used in the stationary 
problem becomes time-dependent, i.e., it is multiplied by 
$\exp(-i \omega n \triangle t)$. 
In case of the second-order discretization, 
$\omega$ is related to the total energy $E$ according to 
\eqref{eq:discrete_E_omega_relation}. Otherwise we simply use $\omega=E/\hbar$.
A time-dependent incoming wave can be included analogously to the one-dimensional 
case outlined in Section \ref{subsec:time_dependent_incoming_waves_1d}.
In particular, we refer to the modifications in 
\eqref{eq:modifications_inc_wave_pml_1d}.
For the sake of brevity, we omit the details.
The storage requirements for the PML are in $O(M J_1^\star)$, where
$J_1^\star$ denotes the number of indices $j_1$ with 
$X_{j_1}^\mathrm{PML}<0$ or $X_{j_1}^\mathrm{PML}>L_1$.
Contrary to DTBC, the computational time required to perform one time step 
is constant.

\begin{figure}
	\centering
  	\includegraphics[width=75mm]{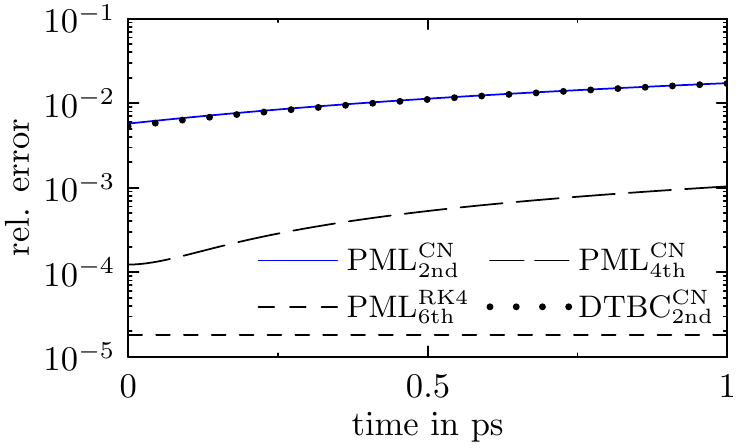}
  	\hspace{2.5mm}
  	\includegraphics[width=75mm]{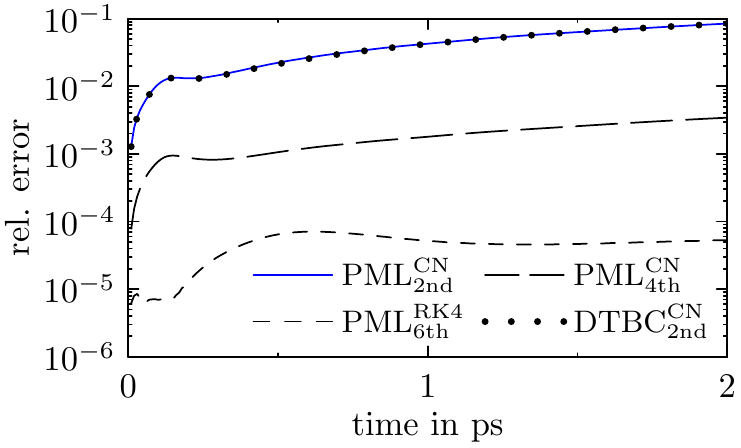}
  	\caption{
  	{\em Left:}
  	Relative errors of the numerical solutions as a function of time 
	in a transient scattering state experiment.
  	{\em Right:}
  	Time evolution of the numerical errors corresponding to two superimposed 
	wave packets in a straight waveguide with parabolic cross section.}
  	\label{fig:rel_errors_transient_scattering_wave_packets_2d}
\end{figure}

\medskip\noindent
{\bf Simulations.}
In the simple example of the straight waveguide \eqref{eq:waveguide_parabolic},
the solution of the transient scattering state problem outlined above is given by
$$
  \psi(x_1,x_2,t) = \exp(i k x_1-i \omega t) \chi^{(0)}(x_2), 
$$
where $\omega = E / \hbar$, $E = \Ekininc + E^{(0)}$, 
$\Ekininc = \hbar^2 k^2 / (2 m^\star)$. Here,
$\chi^{(0)}$ and $E^{(0)}$ are substituted according to 
\eqref{eq:parabolic_cross_sectional_eigenstates}.
We compare the exact solution with the results of the different numerical solvers.
For this, we set $[0,L_1]\times[0,L_2] = [0,120]\times[0,60]\,\mathrm{nm}^2$.
The kinetic energy of the incoming electrons is $\Ekininc=21.5$\,meV.
For the time-integration, we either employ the Crank-Nicolson method or 
the classical Runge-Kutta method.
In the former case, we use $\triangle x = 0.5$\,nm and $\triangle t = 0.25$\,fs.
We note that this choice is rather expensive
since a linear system of equations has to be solved in each time step.
For the Runge-Kutta method, we use the same spatial mesh size but a smaller
time step size of $\triangle t\lesssim 0.05$\,fs is needed to ensure stability.
However, despite of this small time step size, the computation times are the 
shortest compared to those of all other employed solvers.
The relative errors are shown in Figure
\ref{fig:rel_errors_transient_scattering_wave_packets_2d} (left).

The numerical methods derived above can be easily modified to allow 
for simulations of wave packets in quantum waveguides.
For instance, to simulate wave packets using the solver corresponding to DTBC,
we set $c_0^{(m)}=c_1^{(m)}=c_{J_1-1}^{(m)}=c_{J_1}^{(m)}=0$ 
in \eqref{eq:non_hom_dtbc_modes} for all $m=0,\ldots,M-1$.
In case of PML, we simply omit the incoming wave at the device contact.

In the following experiment, we consider the time evolution of
two superimposed wave packets in the straight waveguide with parabolic cross section.
More precisely, we start the numerical calculations using the exact solution
$$
  \psi(x_1,x_2,t) := \xi_p(x_1,t) \exp(-i E^{(0)} t / \hbar) \chi^{(0)}(x_2) 
  + \xi_p(x_1,t) \exp(-i E^{(1)} t / \hbar) \chi^{(1)}(x_2),
$$
at $t=0$\,ps, where $\xi_p$ is defined in \eqref{eq:gaussian_1d}
with the parameters $\sigma = 7.5$\,nm, $x_0 = L_1/2$, and 
$k_p=\sqrt{2 m^\star\,21.5\mathrm{\,meV}}/\hbar$.
The transversal waveguide eigenstates and eigenvalues are defined in 
\eqref{eq:parabolic_cross_sectional_eigenstates}.
We set $[0,L_1]\times[0,L_2] = [0,160]\times[0,60]\,\mathrm{nm}^2$.

The relative errors are depicted in Figure 
\ref{fig:rel_errors_transient_scattering_wave_packets_2d} (right).
We see that $\mathrm{DTBC}_\mathrm{2nd}^\mathrm{CN}$ and 
$\mathrm{PML}_\mathrm{2nd}^\mathrm{CN}$ yield very similar results 
which is in contrast to the one-dimensional results shown in Figure 
\ref{fig:rel_errors_wave_packets_1d}.
In this simulation, the wave packet is confined to a narrow channel
and consequently, the wave function does not only oscillate in the 
longitudinal direction but also in the vertical direction.
For the given cross section of the waveguide  
(characterized by $\omega_\ast = 0.5 \times 10^{14}\mathrm{\,s}^{-1}$),
we find $E^{(0)}\approx 16.45$\,meV and $E^{(1)}\approx 49.36$\,meV.
The oscillation corresponding to the second energy cannot 
be easily resolved by second-order methods at the given spatio-temporal 
resolution. As a consequence, the discretization of the boundary conditions plays a 
minor role and thus $\mathrm{DTBC}_\mathrm{2nd}^\mathrm{CN}$ and 
$\mathrm{PML}_\mathrm{2nd}^\mathrm{CN}$ yield approximately the same results.
We note that for decreasing $\omega_\ast$, the confinement in the $x_2$-direction 
becomes less important.
In such a case, the numerical results would indeed be similar to the numerical 
results of the one-dimensional calculations considered in 
Section \ref{sec:one_dimensional_simulations}.


\subsection{Details of the implementation}

All simulations of this article are realized in the {\em Python} programming 
language using the numerical tools available in {\em SciPy} \cite{JOP01}.
One of the crucial steps is the assembling of the sparse matrices 
related to the solvers using DTBC.
Since these matrices contain small but dense submatrices, one needs to 
proceed carefully. This task can be realized efficiently by using
fast routines to convert dense to sparse matrices 
followed by a series of sparse {\tt vstack} and {\tt hstack} operations.
Another crucial step is the computation of the discrete convolution terms in 
\eqref{eq:non_hom_dtbc_modes_a} and \eqref{eq:non_hom_dtbc_modes_b}.
To this end, we employ parallelized {\em C}-functions which can be included 
easily in {\em Python} programs.
All linear systems are solved using direct solvers for sparse matrices
which is the most time-consuming part in the stationary as well as in the 
transient algorithms.
We note that no linear system needs to be solved if the time evolution is 
computed from the explicit Runge-Kutta time-integration method.
The timings reported below correspond to an 
Intel Core i7-4770K CPU @ 3.50GHz $\times$ 8.


\section{The Aharonov-Bohm effect}\label{sec.AB}

\begin{figure}
  	\begin{centering}
  	\includegraphics[width=150mm]{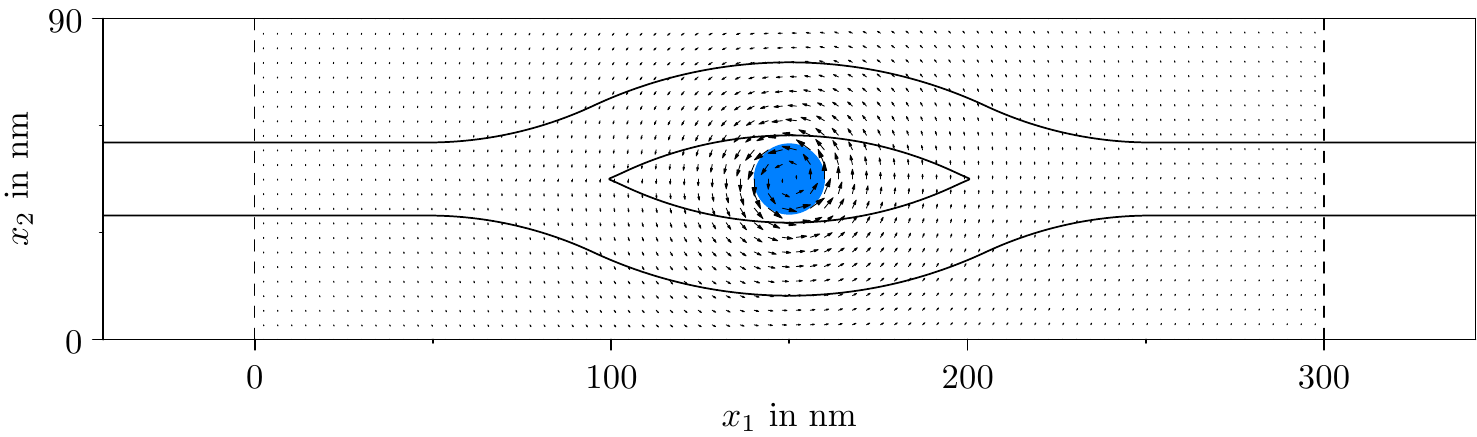}
  	\caption{
  	Vector potential (arrows) and magnetic field (blue disk) in the ring-shaped  
	quantum waveguide device depicted in Figure
  	\ref{fig:potential_energy_aharonov_bohm_ring}.
  	The black solid line indicates an isoline of the potential energy at $200$\,meV.} 
  	\label{fig:vector_potential}
  	\end{centering}
\end{figure}

We consider the ring-shaped quantum waveguide depicted in 
Figure \ref{fig:potential_energy_aharonov_bohm_ring},
but now we apply an additional homogeneous magnetic field, which is perpendicular
to the $(x_1,x_2)$-plane and which is assumed 
to vanish outside a circle of radius $r_0=10$\,nm centered at $(L_1/2,L_2/2)$.
For a charged particle in an electromagnetic field, 
the Hamiltonian of the Schr\"odinger equation reads
\begin{equation}
  \label{eq:hamiltonian_charge_in_em_field}
  \hat{H} = \frac{1}{2 m^\star}(\hat{p} - q A)^2  + q \Phi, \quad 
	\hat{p} = -i \hbar \nabla,
\end{equation}
where $q$ is the particle charge and $A$, $\Phi$ denote the vector and 
scalar potential, respectively.
The electric and magnetic fields are expressed in terms of $A$ and $\Phi$
via $E = - \partial_t A - \nabla \Phi$ and $B = \nabla \times A$.
In the simulations below, the vector potential is defined by
\begin{equation}
  \label{eq:vector_potential}
  A(x_1,x_2,x_2) := B_0 \tilde{A}(x_1-L_1/2,x_2-L_2/2,x_3), 
\end{equation}
with
\begin{equation*}
  \begin{aligned}
  \tilde{A}(x_1,x_2,x_3) &= 
  \begin{cases}
  \frac12 (-x_2,x_1,0)^\top, & \mbox{if } \sqrt{x_1^2 + x_2^2} \leq r_0, \\
  \frac12 r_0^2(x_1^2+x_2^2)^{-1}(-x_2,x_1,0)^\top, 
	& \mbox{if } \sqrt{x_1^2 + x_2^2} > r_0.
  \end{cases}
  \end{aligned}
\end{equation*}
Accordingly, the magnetic field reads as
\begin{equation*}
  B = \begin{cases}
  (0,0,B_0)^\top, & \mbox{if } \sqrt{(x_1-L_1/2)^2 + (x_2-L_2/2)^2} \leq r_0, \\
  (0,0,0)^\top, & \mbox{if } \sqrt{(x_1-L_1/2)^2 + (x_2-L_2/2)^2} > r_0.
  \end{cases}
\end{equation*}
An illustration of $\tilde{A}$ and $B$ is given in Figure \ref{fig:vector_potential}.
The scalar potential $\Phi$ is defined via $V= -e \Phi$. 

In order to fix $A$ and $\Phi$, we impose the Coulomb gauge condition
$\nabla\cdot A=0$. Then, using the electron charge $q=-e$, 
the Hamiltonian in \eqref{eq:hamiltonian_charge_in_em_field} becomes
\begin{equation}
  \label{eq:hamiltonian_electron_in_em_field}
  \hat{H} = -\frac{\hbar^2}{2 m^\star} \Delta 
  - i \frac{e \hbar}{m^\star} A \cdot\nabla 
  + \frac{e^2}{2 m^\star} A^2
  + V.
\end{equation}
Thus, to include the magnetic field into the numerical methods developed in 
Section~\ref{sec:quantum_waveguide_simulations}, 
we need to discretize two additional terms.
The boundary conditions remain unchanged provided that the vector potential 
vanishes in the exterior domains.
In the numerical simulations, we set $A(x_1,x_2,x_3)=0$ for $x_1<\delta$ and 
$x_1>L_1-\delta$ ($\delta\approx 2.5$\,nm) which is not consistent with
\eqref{eq:vector_potential}.
However, the development of DTBC in the presence of the vector potential in 
\eqref{eq:vector_potential}, if at all possible, would be a challenging task.
Here, we simply choose $L_1$ such that $|A(x_1,x_2,x_3)|$ becomes small for 
$x_1<\delta$ and $x_1>L_1-\delta$.
As a result, the modelling error becomes small too, which can be 
verified by further increasing $L_1$.
The term $e^2 A^2/(2m^\star)$ in \eqref{eq:hamiltonian_electron_in_em_field}
is discretized in the same way as the potential energy.
Using the notations corresponding to the case of DTBC described in Section 
\ref{sec:quantum_waveguide_simulations}, the
convection term $-(i e \hbar / m^\star) A\cdot \nabla$ is discretized as follows:
$$
  A \cdot \nabla \approx \diag(A_1) D_{x_1}^{1} \otimes I_{J_2} 
	+ \diag(A_2) I_{J_1} \otimes D_{x_2}^{1}.
$$
Here, $D_{x_1}^{1}$, $D_{x_2}^{1}$ are defined according to 
\eqref{eq:D_1_2nd_4th_6th} and $A_1$, $A_2$ denote the first and second component 
of the vector potential $A$, respectively, i.e.,
$(A_{1,2})_j = (A_{1,2})_{j_1,j_2}$, where $j=j_1 J_2+j_2$, $j_1=0,\ldots,J_1$, 
and $j_2=0,\ldots,J_2$.
The new terms can be easily included into the stationary and transient methods 
outlined in Section~\ref{sec:quantum_waveguide_simulations}.

Let us consider the stationary case first.
In particular, we are interested in the 
transmission probability \eqref{eq:transmission_probability}
as a function of the magnetic flux $\Phi_B = B_0 \pi r_0^2$
(not to be confused with the scalar potential $\Phi$). 
To this end, we compute scattering state solutions of the stationary Schr\"odinger 
equation using the Hamiltonian \eqref{eq:hamiltonian_electron_in_em_field}.
Like in Section~\ref{sec:quantum_waveguide_simulations} we consider electrons injected 
at the left terminal traveling to the right.
While the magnetic flux is increased in each calculation, we keep the 
kinetic energy of the incoming electrons fixed at $\Ekininc=21.5$\,meV.

\begin{figure}
	\centering
  	\includegraphics[width=75mm]{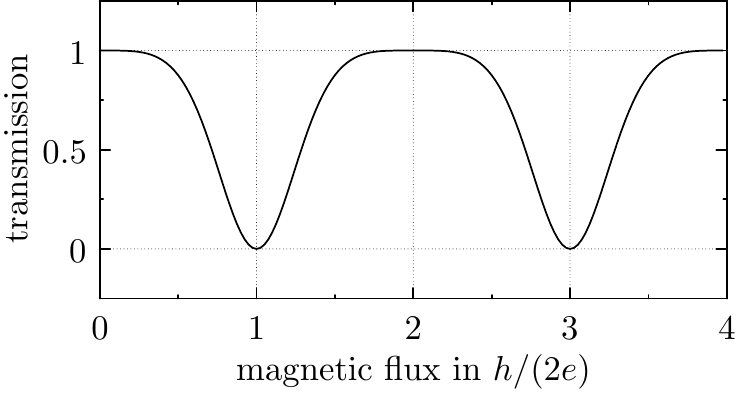}
  	\caption{
	Transmission probability as a function of the magnetic flux in 
	multiples of the flux quantum. 
	The energy of the incoming electrons is fixed at $\Ekininc=21.5$\,meV.
	}
  	\label{fig:transmissions_ab_ring2}
\end{figure}

Figure \ref{fig:transmissions_ab_ring2} shows the computed transmission probability 
as a function of the magnetic flux in multiples of the flux quantum $\Phi_0 = h/(2e)$.
The transmission probability oscillates almost perfectly with increasing values 
of the magnetic flux. These oscillations are a manifestation of the well-known 
Aharonov-Bohm effect \cite{AhBo59}.
In fact, quantum mechanics implies that electrons traveling along a path $P$,
along which $\nabla \times A=0$, accumulate a phase shift
$\varphi = -(e/\hbar) \int_P A \cdot dx$.
Hence, using Stokes' theorem, the phase difference between the beam of electrons 
taking the upper path $P_1$ and the beam taking the lower path $P_2$ 
(see Figure \ref{fig:vector_potential}) is given by
$$
  \varphi_1-\varphi_2 
  = -\frac{e}{\hbar} \left( \int_{P_1} A \cdot dx -\int_{P_2} A \cdot dx \right)
  = -\frac{e}{\hbar} \Phi_B.
$$
Therefore, the interference of the two electron beams depends solely on the 
enclosed magnetic flux. The remarkable fact is that the electrons are affected 
by the vector potential even in regions where the magnetic field is zero 
which is in strong contrast to classical mechanics.

\begin{figure}
	\begin{centering}
	\includegraphics[width=150mm]{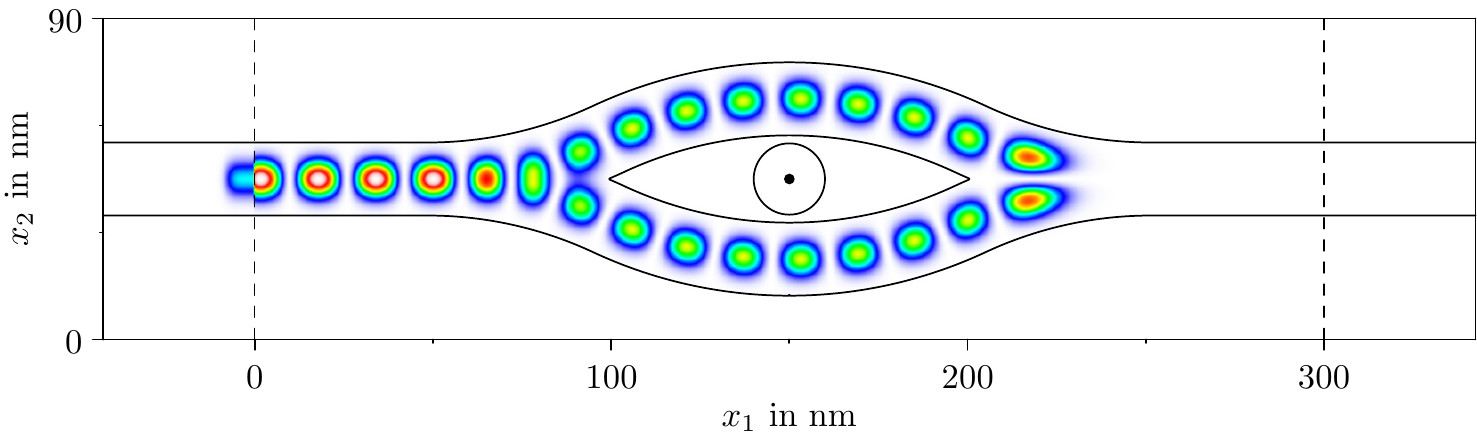}
  	\caption{
  	Scattering state solution in a ring-shaped quantum device using the
  	same colormap as in Figure \ref{fig:scattering_state_ab_ring_E_kin_inc_21_5_meV}.
  	The kinetic energy of the electrons injected at the left contact is
  	$\Ekininc=21.5$\,meV and 
  	the magnetic flux through the encircled area is given by $\Phi_B = h / (2 e)$.
  	The black solid line indicates an isoline of the potential energy at $200$\,meV.}
  	\label{fig:scattering_state_ab_ring_destructive}
  	\end{centering}
\end{figure}

As an example, Figure \ref{fig:scattering_state_ab_ring_destructive} shows the 
scattering state solution for $\Phi_B = h/(2e)$ corresponding to the destructive 
interference condition $|\varphi_1-\varphi_2| = \pi$.
For the calculations, we employed second-order spatial discretizations 
with $\triangle x = 0.5$\,nm
corresponding to linear systems of size
$N_\mathrm{DTBC} = 66\,264$ and $N_\mathrm{PML} = 79\,694$.
The relative difference between the solutions computed using DTBC and PML
is approximately $2\times 10^{-3}$.
We note that this value is almost independent from the magnetic flux.

\begin{figure}[b]
	\centering
  	\includegraphics[width=75mm]{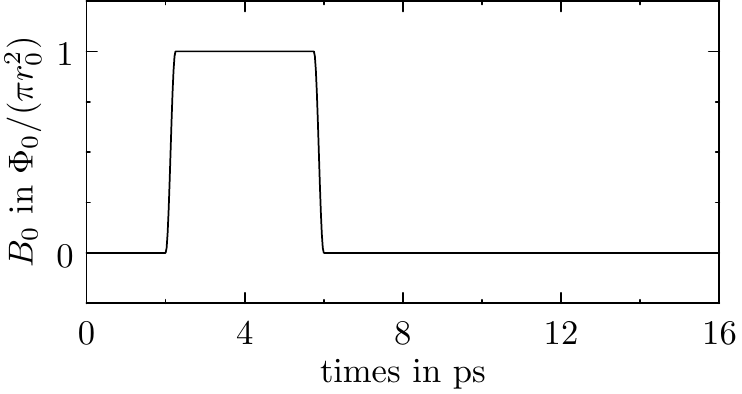}
  	\hspace{2.5mm}
  	\includegraphics[width=75mm]{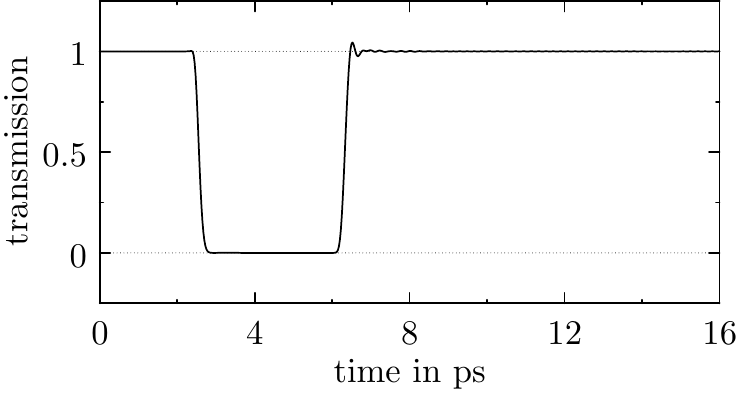}
  	\caption{
  	Magnetic field (left) and transmission probability (right) as a function of time.}
  	\label{fig:B0_and_transmission_of_time}
\end{figure}

\begin{figure}[ht]
  	\centering
  	\includegraphics[width=75mm]{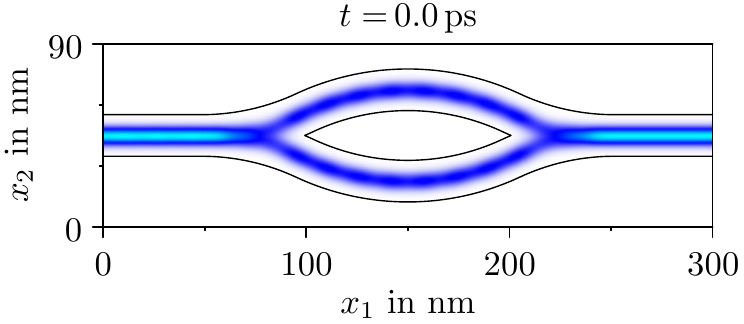}
  	\includegraphics[width=75mm]{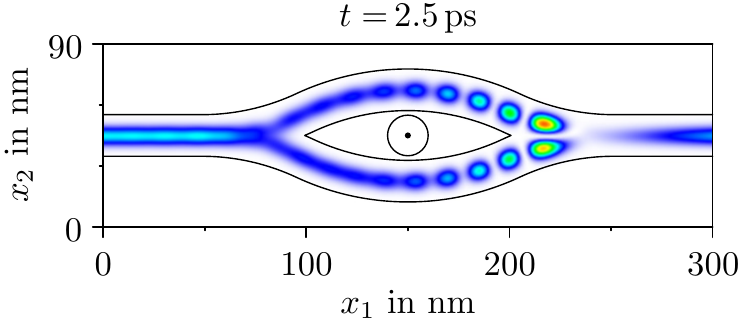} \\
  	\vspace{2.0mm}
  	\includegraphics[width=75mm]{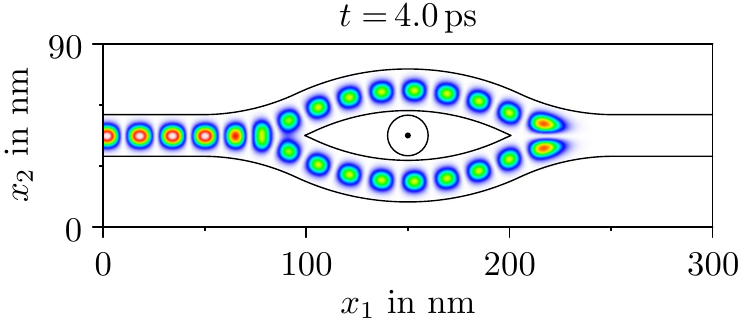}
  	\includegraphics[width=75mm]{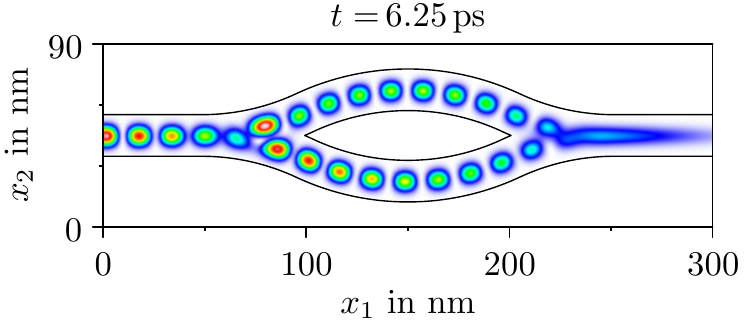} \\
  	\vspace{2.0mm}
  	\includegraphics[width=75mm]{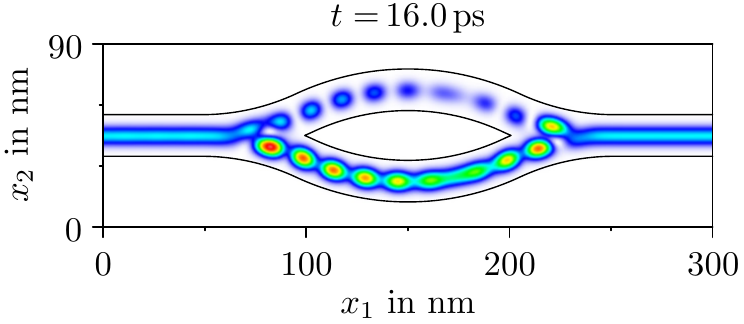}
  	\includegraphics[width=75mm]{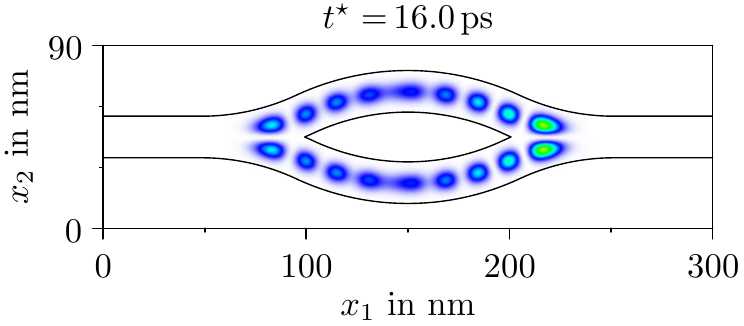}
  	\caption{
	Time evolution of the wave function (probability density) in a
  	transient scattering state experiment with time-dependent magnetic field.
  	The colormap is the same as in Figure 
	\ref{fig:scattering_state_ab_ring_E_kin_inc_21_5_meV}.
	}
	\label{fig:ab_ring_time_evolution}
\end{figure}

We finally turn our attention to the Aharonov-Bohm effect in the transient regime.
In particular, we consider a transient scattering state experiment with
the time-dependent vector potential \eqref{eq:vector_potential}.
More precisely, we let $B_0$ change in time as illustrated in 
Figure \ref{fig:B0_and_transmission_of_time} (left).
Like in the stationary case, we keep the kinetic energy of the 
incoming electrons fixed at $\Ekininc=21.5$\,meV.
Initially, the magnetic field is switched off and thus the initial wave function 
corresponds to the scattering state shown in 
Figure~\ref{fig:scattering_state_ab_ring_E_kin_inc_21_5_meV}.
The time evolution of the wave function is illustrated in 
Figure \ref{fig:ab_ring_time_evolution}.
Figure \ref{fig:B0_and_transmission_of_time} (right) presents the transmission 
probability as a function of time.

Between $t=2$\,ps and $t=2.25$\,ps, the magnetic field is increased to a value of 
$B_0=\Phi_0/(\pi r_0^2)$, which corresponds to the destructive interference 
condition considered in the stationary case.
The effect on the wave function becomes apparent with a short delay 
(Fig.~\ref{fig:ab_ring_time_evolution}, $t=2.5$\,ps) and, as expected, 
the transmission probability decreases.
Another $1.5$\,ps later, the wave function has effectively adopted the scattering 
state corresponding to the destructive interference condition
(Fig.~\ref{fig:ab_ring_time_evolution}, $t=4$\,ps).
During the period from $t=6$\,ps to $t=6.25$\,ps, the magnetic field is turned off.
Again the wave function starts to evolve in a rather wild manner 
(Fig.~\ref{fig:ab_ring_time_evolution}, $t=6.25$\,ps).
Only a short time later, the transmission probability recovers its old value.
However, the wave function does not return to its initial state.
Instead we observe strong oscillations restricted to the interior part 
of the waveguide.
Even after another $9.75$\,ps, these oscillations are still present 
(Fig.~\ref{fig:ab_ring_time_evolution}, $t=16$\,ps) and in fact they will
persist for all time.

\begin{figure}
	\centering
  	\includegraphics[width=75mm]{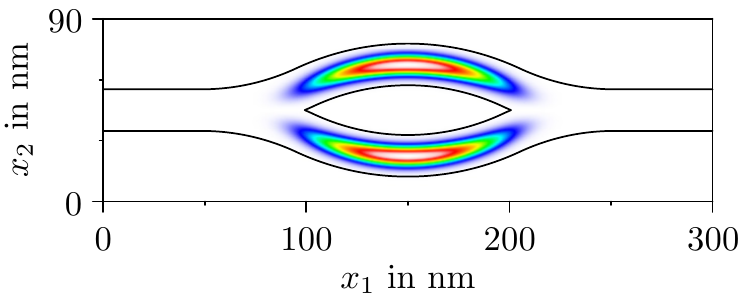}
  	\hspace{0mm}
  	\includegraphics[width=75mm]{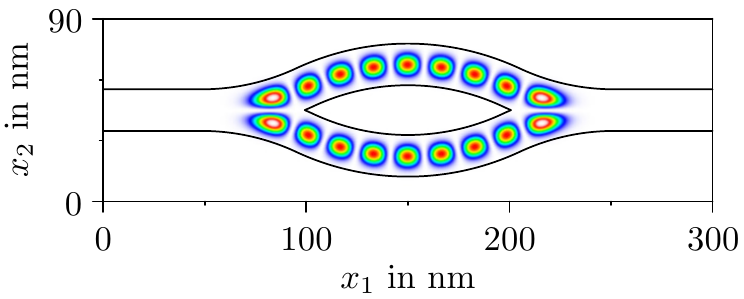}
  	\caption{
  	Eigenfunctions to the stationary Schr\"odinger equation corresponding 
	to the eigenvalues $E\approx 33.2$\,meV and $E\approx 52.8$\,meV.
  	The probability densities are scaled by their maximum values.}
  	\label{fig:bound_states_ab_ring}
\end{figure}

This can be demonstrated by decreasing the amplitude of the incoming wave 
after the second switching operation (which is realized easily in case of PML).
A short time later, we are left with an oscillating wave packet which is bound 
to the region of the ring.
The result of such a numerical simulation is depicted in 
Figure \ref{fig:ab_ring_time_evolution} ($t^\star=16$\,ps).
In this particular example, the amplitude of the incoming wave was multiplied by 
a factor of $0.999$ before each time step with $t\geq 12$\,ps.
The remaining wave packet is supposed to be a superposition 
of several eigenstate solutions to the stationary Schr\"odinger equation. 
We note that the probability density of a superposition of several eigenstates 
is time-dependent.
Two eigenstates corresponding to the eigenvalues $E\approx 33.2$\,meV
and $E\approx 52.8$\,meV are shown in Figure \ref{fig:bound_states_ab_ring} 
but in fact, there are many more eigenstates.
The oscillations superimposed to the final wave function 
are caused by the fast variation of the magnetic field.
When $B_0$ is varied quasi-adiabatically 
(on the time-scale of about $10$\,ps), 
no bound states are excited and hence no oscillations emerge.

The simulation was carried out using the second-order Crank-Nicolson scheme.
Alternatively we used the Runge-Kutta time-integration method in combination 
with a sixth-order spatial discretization and PML.
In case of the Crank-Nicolson scheme, we employed a spatial resolution of 
$\triangle x=1$\,nm and a time-step size of $\triangle t=0.5$\,fs.
Compared to the simulations presented in the previous sections, this discretization 
appears to be rather coarse but it should be noted that
in the example presented above, we have to solve $32\,000$ linear systems of size
$N_\mathrm{DTBC}=16\,536$ and $N_\mathrm{PML}=20\,046$.
If the time evolution is computed using the classical Runge-Kutta
method, no linear systems needs to be solved at all (except for the initial 
scattering state).
However, the resulting numerical method is only conditionally stable and hence
we cannot use the same time step size as for the Crank-Nicolson method.
At a spatial resolution of $\triangle x = 1$\,nm, we are forced to use a time step 
size not larger than $\triangle t=0.25$\,fs which agrees well 
with the condition $\kappa = 9 m^\star/(16 \sqrt{2} \hbar)$
discussed in Section \ref{subsec:wave_packets_1d}.

Comparing the simulations on a reduced mesh 
(see Section \ref{sec:quantum_waveguide_simulations}) with the simulations using 
the full grid, we found that the results are practically the same.
Hence, we did not change the underlying physics by imposing artificial Dirichlet 
boundary conditions at a distance too close to the center of the waveguides.
Eliminating all finite-difference equations corresponding to grid points 
where the wave function is effectively zero, reduces the size of the linear systems 
by more than 40 percent and, as a consequence, the time needed to solve the 
linear systems decreases by more than 50 percent.

\begin{figure}
  	\centering
  	\includegraphics[width=75mm]{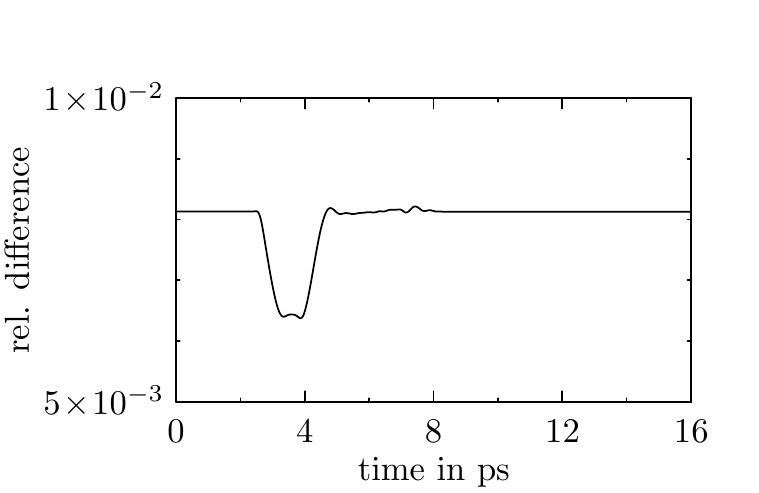}
  	\hspace{2.5mm}
  	\includegraphics[width=75mm]{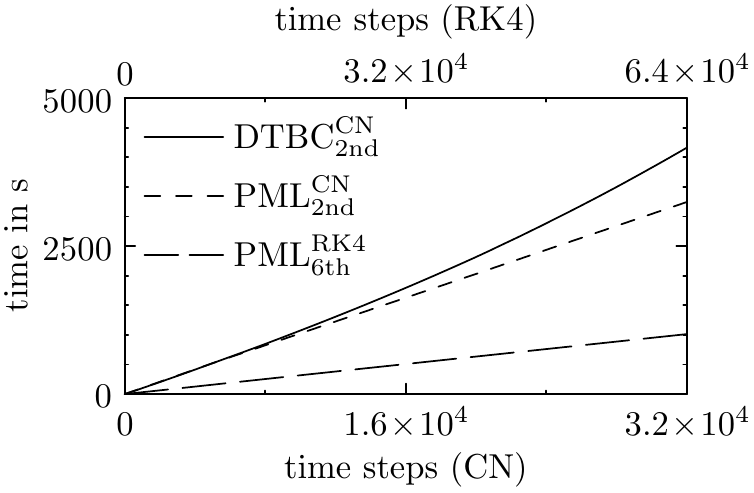}
  	\caption{
  	{\em Left:}
  	Relative difference between the numerical solutions of 
	$\mathrm{DTBC}_\mathrm{2nd}^\mathrm{CN}$ and 
	$\mathrm{PLM}_\mathrm{2nd}^\mathrm{CN}$ as a function of time.
  	{\em Right:}
  	Computing times as a function of the number of time steps.}
	\label{fig:timings_and_rel_differences_ab_ring}
\end{figure}

The relative difference between the numerical solutions obtained by
$\mathrm{DTBC}_\mathrm{2nd}^\mathrm{CN}$ and $\mathrm{PML}_\mathrm{2nd}^\mathrm{CN}$
is shown in Figure \ref{fig:timings_and_rel_differences_ab_ring} (left).
Compared to the relative differences seen in the simulations of the previous 
sections, the difference of both numerical solutions is relatively large.
This is due to the coarse spatial resolution of $\triangle x=1$\,nm which is 
twice as large as the resolution employed before.
Indeed, a spatial resolution of $\triangle x=0.5$\,nm yields the same level of 
accuracy as in the previous sections but the simulation becomes 
extremely time-consuming.

Finally, we report the computing times corresponding to the above mentioned 
methods in Figure~\ref{fig:timings_and_rel_differences_ab_ring} (right).
Initially, $\mathrm{DTBC}_\mathrm{2nd}^\mathrm{CN}$ performs slightly better than
$\mathrm{PML}_\mathrm{2nd}^\mathrm{CN}$.
However, the computing time of $\mathrm{DTBC}_\mathrm{2nd}^\mathrm{CN}$ scales 
quadratically with the number of time steps and thus,
simulations involving even more time steps become very expensive.
We note that this problem could be overcome using approximations of the DTBC 
which can be evaluated more quickly; see, e.g., \cite{AES03}.
As expected, the computing times of the solvers using PML increase linearly.
Interestingly, $\mathrm{PML}_\mathrm{6th}^\mathrm{RK4}$ performs significantly 
faster than the other methods.
Moreover, it is the most accurate method
and its implementation is relatively simple.

A movie showing the temporal evolution of the wave function corresponding to the 
last numerical experiment is available at
\url{http://www.asc.tuwien.ac.at/~juengel}.

\section*{Acknowledgements}
The authors acknowledge partial support from   
the Austrian Science Fund (FWF), grants P20214, P22108, I395, and W1245.


\end{document}